\documentclass[11pt,reqno]{amsart}

\usepackage{etoolbox}

\newbool{bForSubmission}
\booltrue{bForSubmission}

\newbool{bExperimental}

\newbool{bForUs}

\usepackage{amscd,amsmath}
\usepackage{mathrsfs}
\usepackage{amsfonts}

\usepackage{amssymb, bm, xspace}
\usepackage{enumerate}

\usepackage{enumitem}

\usepackage{setspace}

\usepackage[OT2,OT1]{fontenc}

\usepackage{datetime}

\usepackage[dvipsnames]{xcolor}

\usepackage[margin=1.2in]{geometry}

\usepackage{mathtools}

\mathtoolsset{centercolon}

\usepackage{stackengine}
\stackMath

\usepackage{tcolorbox}
\usepackage{tikz}

\newbool{HaveBBM}
\booltrue{HaveBBM}

\ifbool{HaveBBM}{
	\usepackage{bbm}
    }
    {
    }

\usepackage[pagebackref=true, colorlinks=true, citecolor=blue]{hyperref}

\usepackage{cleveref}

\tcbuselibrary{skins}
\usetikzlibrary{shadings}
\tcbset{
    myimage/.style={
        enhanced,
        overlay={
            \begin{scope}[shift={([xshift=1mm, yshift=7mm]frame.north west)}]
            \end{scope}}}}

\tcbset{
    skin=enhanced,
    fonttitle=\bfseries,
    interior style={white},
    segmentation style={black,solid,opacity=0.2,line width=1pt}}

\newtcolorbox{TitledBox}[2][]{
    myimage,              
    coltitle=black,       
    colbacktitle=white,   
    title=My title,
    attach boxed title to top center={
        yshift=-3mm,
        yshifttext=-1mm},
    attach boxed title to top left={
        xshift=1cm,
        yshift=-2mm},
    boxed title style={
        size=small},
    title={#2},#1}

\newcommand{\IncludeFigure}[2]{
    \begin{center}
            \includegraphics[scale=#1]{#2.jpg}
    \end{center}
}

\newcommand{\TwoColumn}[6]%
{%
\begin{minipage}[#3]{#1\textwidth}%
#5%
\end{minipage}%
\begin{minipage}[#4]{#2\textwidth}%
#6%
\end{minipage}%
}

\newcommand{\TwoColumnTop}[4]%
{%
\begin{minipage}[t]{#1\textwidth}
#3
\end{minipage}%
\begin{minipage}[t]{#2\textwidth}
#4
\end{minipage}%
}

\newcommand{\PullMarginsIn}[1]%
{%
\begin{minipage}[c]{0.3\textwidth}%
\phantom{x}%
\end{minipage}%
\begin{minipage}[c]{0.4\textwidth}%
#1%
\end{minipage}%
\begin{minipage}[c]{0.3\textwidth}%
\phantom{x}%
\end{minipage}%
}

\newcommand{\Comment}[1]{{\color{Brown}#1}}
\newcommand{\OptionalDetails}[1]{
    \ifbool{bForSubmission}
        {%
        }%
        {\begin{quote}\Comment{\footnotesize
        \medskip
        
        \noindent\textbf{Details not for submission}: \\
        \noindent#1}
        \end{quote}
        }
    }

\newbool{arXivFormat}
\boolfalse{arXivFormat}
    
\newcommand{\IfarXivElse}[2]{
    \ifbool{arXivFormat}
        {#1}{#2}
    }

\renewcommand{\mathbf}[1]{\bm{#1} \textbf{ *** Use bm instead of mathbf ***}}

\newcommand{\eqn}{\begin{eqnarray}}
\newcommand{\een}{\end{eqnarray}}

\newtheorem{theorem}{Theorem}[section]
\newtheorem*{theorem*}{Theorem}				
\newtheorem{prop}[theorem]{Proposition}
\newtheorem{lemma}[theorem]{Lemma}
\newtheorem{cor}[theorem]{Corollary}

\newtheorem{definition}[theorem]{Definition}
\newtheorem{remark}[theorem]{Remark}
\newtheorem*{remark*}{Remark}

\newtheorem{assumption}[theorem]{Assumption}
\numberwithin{equation}{section}

\newcommand{\innp}[1]{\ensuremath{\langle #1 \rangle}}

\newcommand{\BoldTau}
    {{\mbox{\boldmath $\tau$}}}

\newcommand{\bgamma}
    {{\mbox{\boldmath $\gamma$}}}

\newcommand{\BB}[1]{\ensuremath{\mathbb{#1}}}
\newcommand{\R}{{\ensuremath{\BB{R}}}}

\newcommand{\Z}{\ensuremath{\BB{Z}}}

\providecommand{\C}{\ensuremath{\BB{C}}}
\renewcommand{\C}{\ensuremath{\BB{C}}}

\newcommand{\iny}{\ensuremath{\infty}}
\newcommand{\grad}{\ensuremath{\nabla}}
\newcommand{\bigvec}{\overrightarrow}
\newcommand{\vecinv}{\overleftarrow}

\newcommand{\CharFunc}{
    \ifbool{HaveBBM}{
        {\ensuremath{\mathbbm{1}}}
        }
        {
        {\ensuremath{\bm{1}}}
        }
    }

\DeclareMathOperator{\dv}{div} %
\DeclareMathOperator{\curl}{curl} %

\DeclareMathOperator{\image}{image} %
\DeclareMathOperator{\RE}{Re} %
\DeclareMathOperator{\IM}{Im} %
\DeclareMathOperator{\PV}{p.v.} %
\DeclareMathOperator{\BUC}{BUC} %
\newcommand{\prt}{\ensuremath{\partial}}
\newcommand{\brac}[1]{\ensuremath{\left[ #1 \right]}}

\newcommand{\pr}[1]{\ensuremath{\left( #1 \right)}}

\newcommand{\set}[1]{\ensuremath{\left\{ #1 \right\}}}

\DeclarePairedDelimiterX{\norm}[1]{\lVert}{\rVert}{#1}

\DeclarePairedDelimiterX{\abs}[1]{\lvert}{\rvert}{#1}

\newcommand\tenq[2][1]{%
	\def\useanchorwidth{T}%
	\ifnum#1>1%
		\stackunder[0pt]{\tenq[\numexpr#1-1\relax]{#2}}{\scriptscriptstyle\sim}%
	\else%
		\stackunder[1pt]{#2}{\scriptscriptstyle\sim}%
	\fi%
	}

\newcommand{\n}{{\bm{n}}}

\newcommand{\wh}{\widehat}

\renewcommand{\epsilon}{\varepsilon}
\newcommand{\eps}{\ensuremath{\varepsilon}}

\newcommand{\Cal}[1]{\ensuremath{\mathcal{#1}}}
\newcommand{\al}{\ensuremath{\alpha}}

\newcommand{\pdx}[2]{\frac{\prt #1}{\prt #2}}

\newcommand{\diff}[2]{\frac{ d#1}{d#2}}

\newcommand{\ol}{\overline}

\renewcommand{\le}{\leqslant}
\renewcommand{\ge}{\geqslant}

\newcommand{\Ignore}[1]{}

\newcommand{\Obsolete}[1]{}

\newcommand{\Experimental}[1]%
{%
\ifbool{bExperimental}%
	{\bigskip
	\noindent
	\textbf{\color{Brown}*** Start experimental}\\
	#1
}
{
}
}

\definecolor{Correction}{named}{red}

\crefname{cor}{Corollary}{Corollaries} 
									   
\crefname{lemma}{Lemma}{Lemmas}	       

\crefname{section}{Section}{Sections}
\Crefname{section}{Section}{Sections}

\crefname{appendix}{Appendix}{Appendices}
\Crefname{appendix}{Appendix}{Appendices}

\crefname{theorem}{Theorem}{Theorems}
\Crefname{theorem}{Theorem}{Theorems}

\crefname{prop}{Proposition}{Propositions}
\Crefname{prop}{Proposition}{Propositions}

\crefname{conj}{Conjecture}{Conjectures}
\Crefname{conj}{Conjecture}{Conjectures}

\crefname{definition}{Definition}{Definitions}
\Crefname{definition}{Definition}{Definitions}

\crefname{remark}{Remark}{Remarks}
\Crefname{remark}{Remark}{Remarks}

\crefname{figure}{figure}{figures}
\Crefname{figure}{Figure}{Figures}

\crefname{assumption}{Assumption}{Assumptions}
\Crefname{assumption}{Assumption}{Assumptions}

\crefformat{equation}{(#2#1#3)}
\crefrangeformat{equation}{(#3#1#4) through (#5#2#6)}
\crefmultiformat{equation}
    {(#2#1#3)}%
    { and~(#2#1#3)}
    {, (#2#1#3)}
    { and~(#2#1#3)}

\newcommand{\e}{\bm{\mathrm{e}}}

\newcommand{\uu}{{\bm{\mathrm{u}}}}   

\newcommand{\vv}{{\bm{\mathrm{v}}}}

\newcommand{\ww}{{\bm{\mathrm{w}}}}
\newcommand{\x}{{\bm{\mathrm{x}}}}
\newcommand{\y}{{\bm{\mathrm{y}}}}

\newcommand{\UU}{{\bm{\mathrm{U}}}}

\newcommand{\wv}{\widetilde{\vv}}

\newcommand{\wq}{\widetilde{q}}

\newcommand{\PiR}{{\Pi_p}}

\newcommand{\Per}{\Cal{P}er}
\newcommand{\Rep}{\Cal{R}ep}
\newcommand{\KSym}{K_{sym}}

\newcommand{\stardot}{\mathop{* \cdot}}

\DeclareMathOperator{\Res}{Res}

\newcommand{\Cover}{p}
\newcommand{\KInfFunc}{\rho}

\def\Xint#1{\mathchoice
{\XXint\displaystyle\textstyle{#1}}%
{\XXint\textstyle\scriptstyle{#1}}%
{\XXint\scriptstyle\scriptscriptstyle{#1}}%
{\XXint\scriptscriptstyle\scriptscriptstyle{#1}}%
\!\int}
\def\XXint#1#2#3{{\setbox0=\hbox{$#1{#2#3}{\int}$ }
\vcenter{\hbox{$#2#3$ }}\kern-.6\wd0}}

\def\complexint{\Xint{\C}}

\makeatletter
\ifcsname phantomsection\endcsname
    \newcommand*{\qrr@gobblenexttocentry}[5]{}
\else
    \newcommand*{\qrr@gobblenexttocentry}[4]{}
\fi
\newcommand*{\addsubsection}{%
    \addtocontents{toc}{\protect\qrr@gobblenexttocentry}%
    \subsection}
\makeatother

\allowdisplaybreaks

\begin{document}
\newdateformat{mydate}{\THEDAY~\monthname~\THEYEAR}

\title
	[Periodic vortex patches and layers]
	{Contour dynamics and global regularity for periodic vortex patches and layers}

\author[D. Ambrose, F. Hadadifard, and J. Kelliher]
{David M. Ambrose$^{1}$, Fazel Hadadifard$^{2}$, and James P. Kelliher$^{2}$}
\address{$^1$ Department of Mathematics, Drexel University, 3141 Chestnut Street, Philadelphia, PA 19104}
\address{$^2$ Department of Mathematics, University of California, Riverside, 900 University Ave., Riverside, CA 92521}
\email{ambrose@math.drexel.edu}
\email{fazelh@ucr.edu}
\email{kelliher@math.ucr.edu}

\begin{abstract}
We study vortex patches for the 2D incompressible Euler equations.  Prior works on this problem
take the support of the vorticity (i.e., the vortex patch) to be a bounded region.  We instead consider
the horizontally periodic setting.  This includes both the case of a periodic array of bounded vortex patches and the case of vertically bounded vortex layers.  We develop the contour dynamics equation for the
boundary of the patch in this horizontally periodic setting, and demonstrate global $C^{1,\epsilon}$ regularity
of this patch boundary.  In the process of formulating the problem, we consider different notions of periodic 
solutions of the 2D incompressible Euler equations, and demonstrate equivalence of these.
\end{abstract}

\maketitle

\markleft{D.M. AMBROSE, F. HADADIFARD, and J. KELLIHER}

\vspace{-2.5em}

\begin{center}
\medskip
Compiled on {\dayofweekname{\day}{\month}{\year} \mydate\today} at \currenttime
		
\end{center}

\bigskip

\vspace{-2.0em}

\setcounter{tocdepth}{1}
{
\renewcommand\contentsname{}	

\tableofcontents
}

%
%
\section{Introduction}\label{S:Introduction}

\noindent A 2D vortex patch is a solution to the 2D Euler equations for which the vorticity is a constant multiplied by the characteristic function of a domain. We investigate the behavior of vortex patches in an infinite strip periodic in one direction, topologically $S^1 \times \R$, and the corresponding behavior of the vortex patch or layer in the full plane.  Our main results are the extension of 
the $C^{1,\eps}$ global regularity theory for the boundary of the vortex patch to this case, developing and using the 
appropriate contour dynamics equation for this purpose. Here, and throughout, we fix $\eps \in (0, 1)$.

\subsection{The Euler equations}

We can write the 2D incompressible Euler equations (without forcing) on a domain $U$ in vorticity form as
\begin{align}\label{e:EGen}
	\begin{cases}
		\prt_t \omega + \uu \cdot \grad \omega = 0
			&\text{in } \R \times U, \\
		\uu = K[\omega]
			&\text{in } \R \times U, \\
		\omega(0) = \omega^0
			&\text{in } U.
	\end{cases}
\end{align}
Here, $\omega$ is the vorticity---the scalar curl of the velocity field $\uu$. The vorticity is transported by the velocity field as in \cref{e:EGen}$_1$, and the velocity field is recovered from the vorticity field by the constitutive law in \cref{e:EGen}$_2$ so as to be divergence-free and to satisfy any boundary conditions, decay at infinity, or periodicity that might be demanded based, in part, upon the nature of the domain $U$.

Classically, if $U = \R^2$ and the solution has sufficient decay, one uses the Biot-Savart law as the constitutive law:
\begin{align}\label{e:BSR2}
	K[\omega] := K * \omega, \quad
	K(\x)
		:= \grad^\perp \brac{\frac{1}{2 \pi} \log \abs{\x}}
		= \frac{1}{2 \pi} \frac{\x^\perp}{\abs{\x}^2}.
\end{align}
Here, $K$ is the Biot-Savart kernel, which we note lies in $L^1_{loc}(\R^2)$, though $K \notin L^p(\R^2)$ for any $p \in [1, \iny]$. To handle solutions having insufficient spatial decay of the vorticity, we must either find an appropriate substitute for the Biot-Savart law or avoid it entirely by using a velocity, pressure formulation.

\subsection{The plane and the cylinder}

In this paper, we will consider two domains: $U = \R^2$ and $U = \Pi$, the infinite flat periodic strip, $S^1 \times \R \cong \R^2 / \Z \cong \C / \Z$, which we will most often treat in the form
\begin{align}\label{e:Pi}
	\Pi := \brac{-\tfrac{1}{2}, \tfrac{1}{2}} \times \R \text{ with }
			\set{-\tfrac{1}{2}} \times \R
			\text{ identified with }
			\set{\tfrac{1}{2}} \times \R.
\end{align}
We will also find use for the same set as a subset of $\R^2$ or $\C$ without identifying its sides:
\begin{align}\label{e:PiR}
	\PiR := \pr{-\tfrac{1}{2}, \tfrac{1}{2}} \times \R \subseteq \R^2.
\end{align}

\begin{figure}[ht]
\IncludeFigure{0.14}{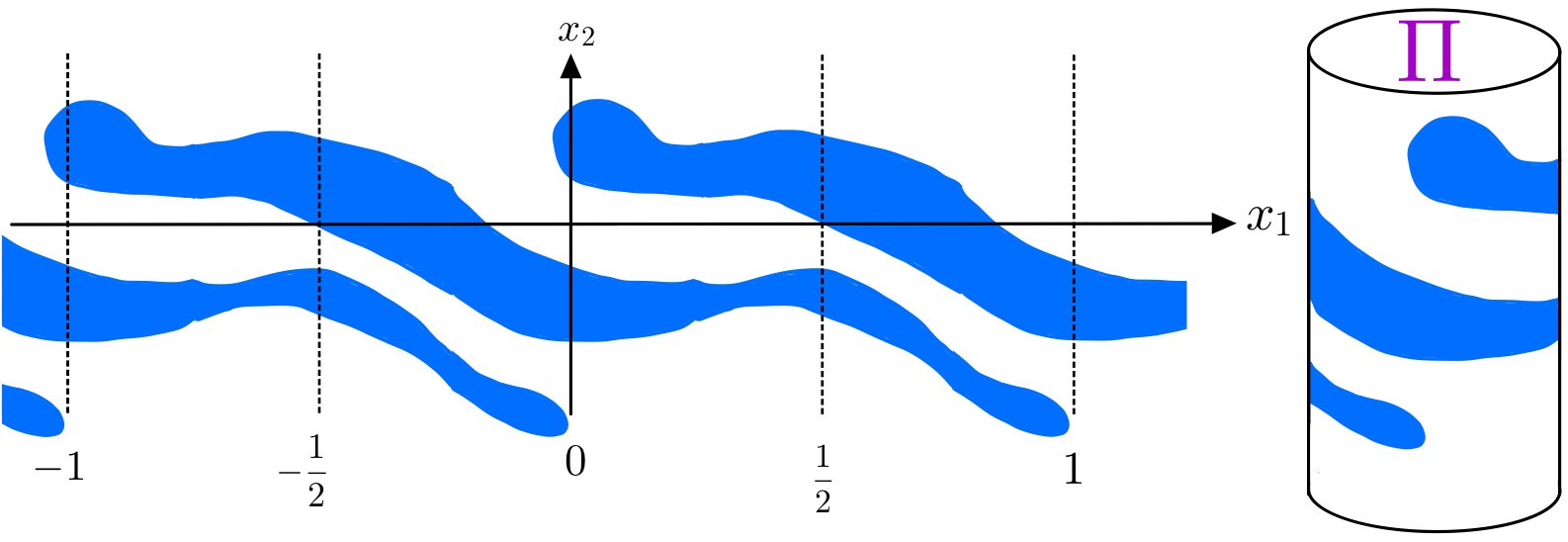}
		
\vspace{-1.2em}
	
\caption{\textit{Example of a periodic vortex patch in $\R^2$ and in $\Pi$}}\label{f:VortexPatch}
\end{figure}

\begin{figure}[ht]

\IncludeFigure{0.14}{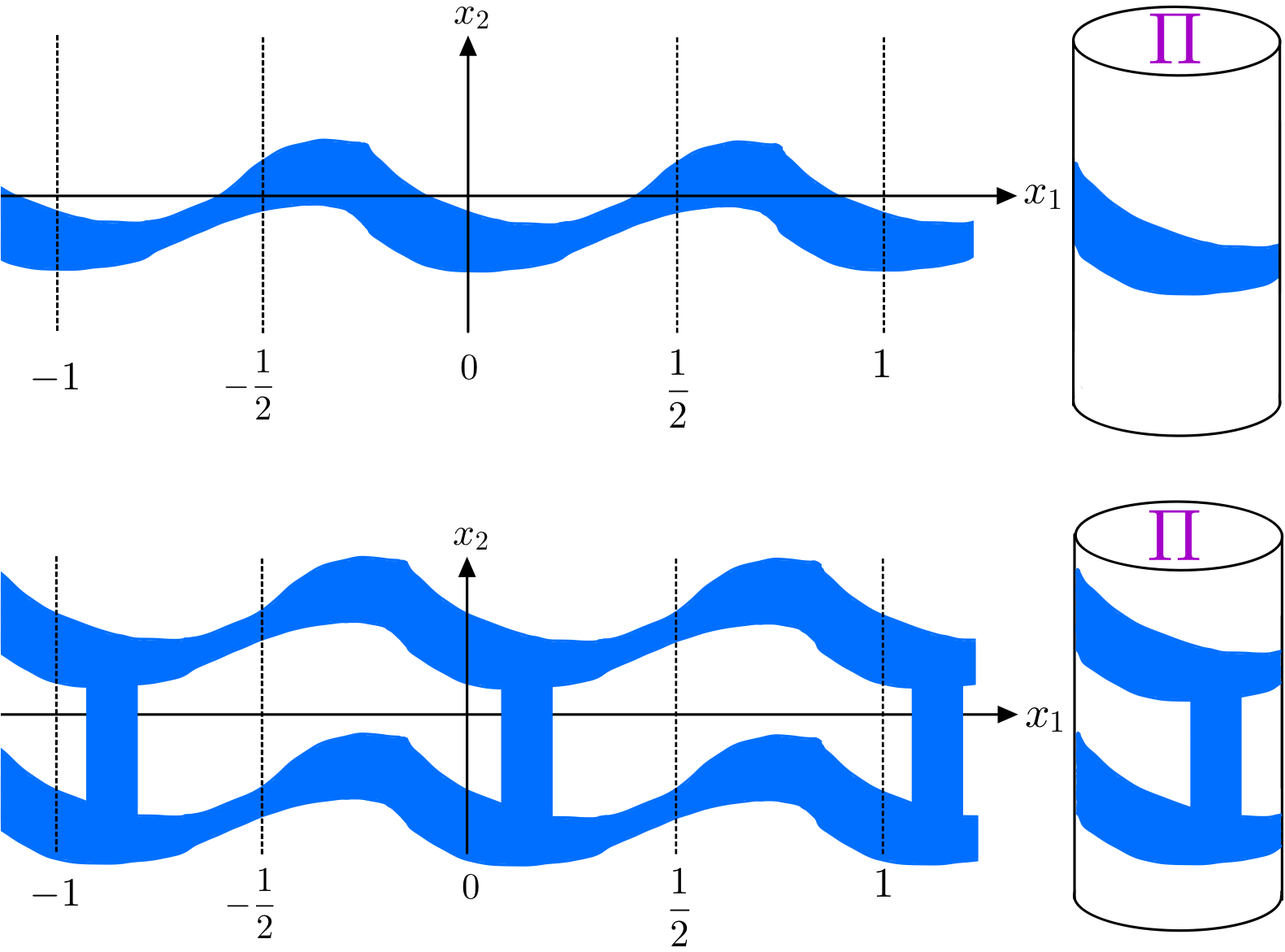}	
\vspace{-1.5em}
	
\caption{\textit{Two examples of a periodic vortex layer in $\R^2$ and in $\Pi$}}\label{f:VortexLayer}
\end{figure}

Suppose we have an initial vorticity $\omega^0 = \CharFunc_\Omega$ for $\Omega$ a bounded domain in $\Pi$. We can periodize it to obtain an initial vorticity in $\R^2$ that is periodic in $x_1$. What results may consist of an infinite number of disconnected domains repeated periodically, one connected, $x_1$-periodic domain, or a combination of each. \Cref{f:VortexPatch} displays an example of a simply connected bounded domain in $\Pi$ yielding an infinite number of copies of the domain in $\R^2$. \Cref{f:VortexLayer} displays two examples of a non-simply connected domain in $\Pi$ producing one domain in $\R^2$ periodically repeating in $x_1$, a so-called \textit{vortex layer}.

On the other hand, we could instead formulate the problem by starting
with an initial vortex patch in $\R^2$ and periodize it in $x_1$. If we can translate the evolution of the patch in $\R^2$ to the evolution in $\Pi$ and back, we can use an understanding of patch behavior in $\Pi$ to gain an understanding of the periodic behavior in $\R^2$. The translation back and forth between $\Pi$ and $\R^2$ is best understood in the more general setting of weak solutions to the 2D Euler equations for bounded vorticity, which includes vortex patch data as a special case.

\subsection{Three types of solutions}

Toward this end, we consider three types of solution to the 2D Euler equations. We summarize the three types of solution briefly now, giving more complete descriptions in later sections.
\begin{itemize}[leftmargin=5em,itemsep=0.5em]
	\item[\textbf{Type 1}]
		Assume that $\uu^0 \in L^\iny(\R^2)$ is divergence-free with
		$\omega^0 := \curl \uu^0 \in L^\iny(\R^2)$ as well.
		Obtain a bounded vorticity, bounded velocity solution to the
		the 2D Euler equations on all of $\R^2$ having initial velocity $\uu^0$ as done by Serfati
		in \cite{Serfati1995Bounded}.
				
	\item[\textbf{Type 2}]
		Assume $\uu^0 \in L^\iny(\Pi)$ is divergence-free with
		$\omega^0 := \curl \uu^0 \in L^\iny(\Pi)$ as well.
		Solve the 2D Euler equations in $\Pi$,
		as done in \cite{AfendikovMielke2005,GallaySlijepcevic2014,GallaySlijepcevic2015}.

	\item[\textbf{Type 3}]
		Let $\omega^0 \in L^\iny(\R^2)$ be compactly supported.
		Solve the 2D Euler equations in vorticity form in all of $\R^2$ with initial vorticity
		$\omega^0$, but  recovering the velocity
		by applying the Biot-Savart law symmetrically to pairs of the periodically
		extended copies of $\omega$. This leads to a replacement Biot-Savart kernel,
		$K_\iny$. 
\end{itemize}

Type 1 and Type 2 solutions are for (potentially) non-decaying velocity and vorticity, but for Type 3 we restrict our attention to vertically decaying solutions, since our primary application is to vortex patch data. Moreover, the convolution $K_\iny * \omega$ cannot be easily defined without some decay assumption.

We will find that all three types of solution are equivalent for a large class of  initial data. Since our primary interest is in vortex patches and layers, we will keep things simple by assuming compact support in $\Pi$. Assuming, then, that $g \in L^\iny_c(\R^2)$---the space of essentially bounded functions with compact support---we define $\Per(g)$ on $\Pi$ by
\begin{align*}
	\Per(g)(\x) = \sum_{n \in \Z} g(\x - (n, 0)),
\end{align*}
noting that for each $\x$ the sum has only finitely many nonzero terms. For any measurable function $f$ on $\Pi$ we define $\Rep(f)$ on $\R^2$ by
\begin{align*}
	\Rep(f)(\x)
		:= f(x_1 - \lfloor x_1 + \tfrac{1}{2} \rfloor, x_2).
\end{align*}	

\begin{figure}[ht]

\IncludeFigure{0.12}{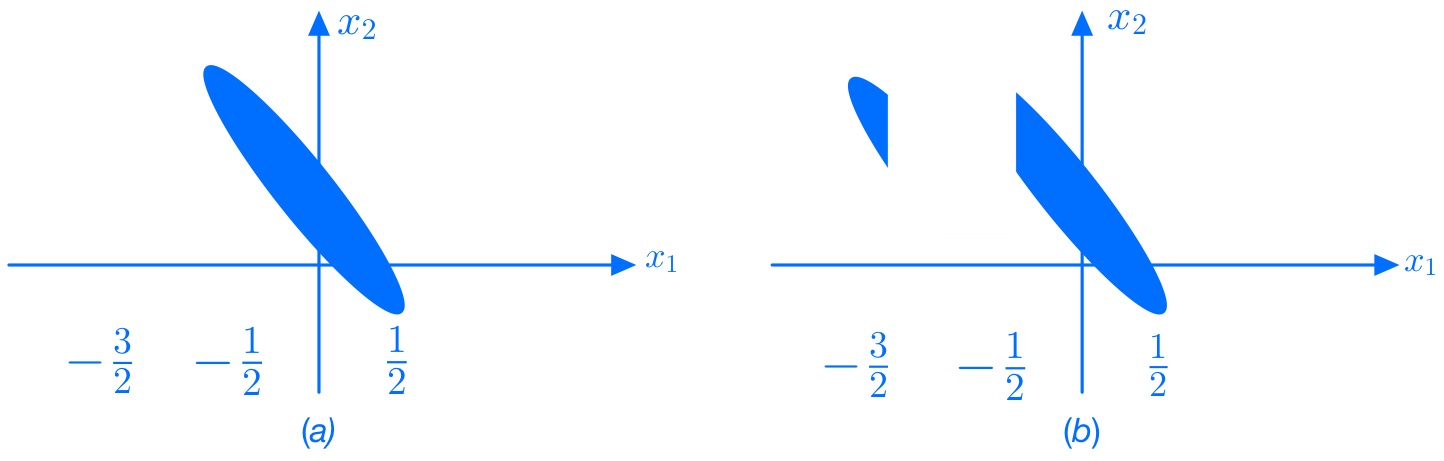}

\vspace{-1.2em}
	
\caption{\textit{Support of two $L^\iny_c(\R^2)$ functions in the same equivalence class}}\label{f:PropertyP}
\end{figure}

\begin{definition}\label{D:Equivalence}
	Two functions $g_1, g_2 \in L^\iny_c(\R^2)$ are equivalent,
	$g_1 \sim g_2$, if $\Per(g_1) = \Per(g_2)$.
	\Cref{f:PropertyP} depicts the support of two functions in the same equivalence class.
\end{definition}

Suppose that $g \in L^\iny_c(\R^2)$, and for purposes of illustration, let us treat it as the characteristic function of a bounded domain (our primary application), whose support is depicted as in either (a) or (b) of \Cref{f:PropertyP}. Below, we construct an initial vorticity from $g$ and depict the support of $\omega^0$ for each type of solution (the time-evolved vorticity being of a similar nature).

\begin{itemize}[leftmargin=5em,itemsep=0.5em]
	\item[\textbf{Type 1}]
		Let $\omega^0 = \Rep(\Per(g))$.
		
		\IncludeFigure{0.12}{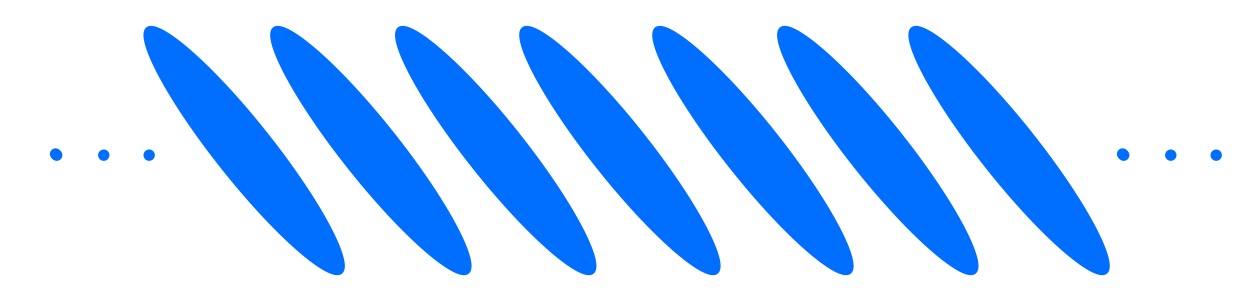}
				
	\item[\textbf{Type 2}]
		Let Let $\omega^0 = \Per(g)$.

		\vspace{-0.5em}
		\IncludeFigure{0.12}{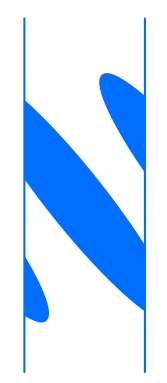}
		
	\item[\textbf{Type 3}]
		Let $\omega^0 = g$.
		The vorticity $\omega$ is transported
		by the flow from the single
		copy of $g$, and so is no longer the curl of
		$\uu$. There are, in effect, multiple
		phantom copies of $g$ matching those of Type 1.
			
		\IncludeFigure{0.12}{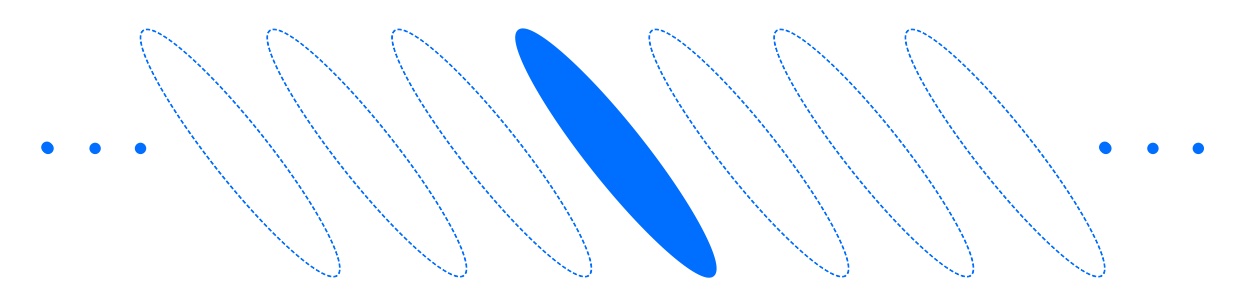}
		\begin{center}
			or
		\end{center}
		\IncludeFigure{0.12}{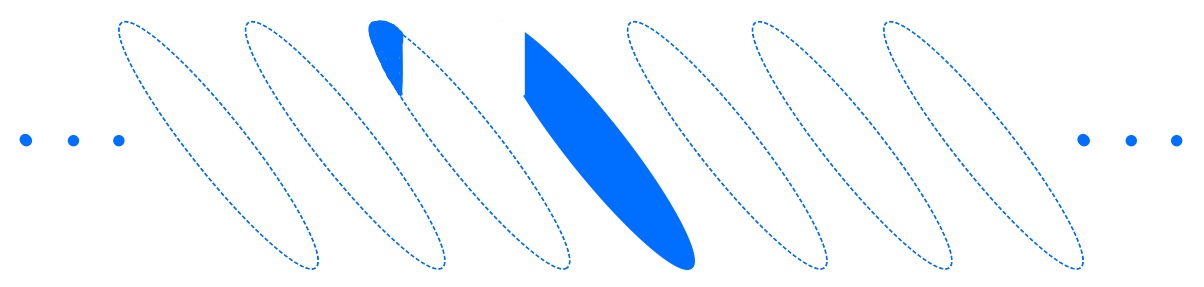}
\end{itemize}

The vorticity $\omega^0$ for Type 1 and 2 do not depend upon the representative for the equivalence class, though Type 3 does. We will find, nonetheless, that the velocity field for solutions of Type 3 is independent of the representative.

It is mentioned in \cite{Gallay2017} that a Type 2 solution is equivalent to a Type 1 solution with periodic velocity and pressure. Following up on this comment, we will show that all three types of solution are equivalent. The equivalence of Type 1 and Type 2 solutions, which applies to a larger class of initial data than we have so far discussed, will rely upon the properties of the pressure required for uniqueness for those two types of solution. The equivalence of Type 3 and Type 2 (and so of Type 1) will rest primarily on showing that solutions of Type 2 reduce to those of Type 3 when the vorticity has sufficient vertical decay. A side benefit of this approach is that it will give the well-posedness of Type 3 solutions. Such a well-posedness result could be obtained by adapting in a fairly straightforward way the approach Marchioro and Pulvirenti take in \cite{MP1984,MP1994} for the 2D Euler equations, except for subtle points regarding the periodicity of the pressure. It is thus more efficient to leverage the technology developed in \cite{AfendikovMielke2005,GallaySlijepcevic2014,GallaySlijepcevic2015}, though it is more than is strictly needed to develop Type 3 solutions alone.

Specializing to vortex patch data, we will then show how the contour dynamics equation (CDE) is adapted from the classical form, which allows the propagation of regularity of the boundary of a vortex patch to be proved, adapting the argument of Bertozzi and Constantin in \cite{ConstantinBertozzi1993}.

\subsection{Prior work}\label{S:ToDo}
Bounded vortex patches evolving under the two-dimensional Euler equations have been well-studied, with global 
regularity of the boundary being established by Chemin \cite{cheminCRAS} and by Bertozzi and Constantin 
\cite{ConstantinBertozzi1993}.  Regularity of the vortex patch boundary can also be seen to follow from a more general
approach studying level sets of the vorticity, establishing striated regularity, as in the work of Chemin \cite{cheminENS}
and Serfati \cite{serfatiStratifee}.  Regularity of bounded vortex patches and/or striated regularity have been established for solutions
of related evolution equations as well, such as aggregation equations \cite{bertozziEtAl}, active transport equations 
\cite{baeKelliher}, and the surface quasi-geostrophic equation and related systems \cite{chaeEtal}, \cite{gancedo},
\cite{kiselev}.  None of these problems consider unbounded vortex patches as in the present work.

There are seemingly fewer papers on the evolution of vortex layers.  
An equation similar to our version of the contour dynamics equation for the motion of the patch/layer boundary was developed in
\cite{pullin}, and was subsequently used in \cite{golubeva} for the study of complex singularities in vortex layers.
(We mention that the version of the contour dynamics equation developed in the present work lends itself to the 
study of global regularity.)
Atassi, Bernoff, and Lichter study the interaction of a point vortex with a vortex layer \cite{bernoff}.  Crowdy gives some exact
solutions of vortex layers interacting with solid boundaries \cite{crowdy}.  Benedetto and Pulvirenti have shown that vortex
layers rigorously approximate vortex sheets in analytic function spaces \cite{benedettoPulvirenti}.
Caflisch, Sammartino, and collaborators have considered vortex layers which are not sharp fronts in a series of papers \cite{caflisch2}, \cite{caflisch1}, \cite{caflisch3}, 
considering how such flows behave in the zero viscosity limit and how such flows may approximate vortex sheets, which represent
a more singular vorticity configuration.  In these works, they take the vorticity to be exponentially decaying (in the vertical direction) away from a core region,
rather than being an indicator function as in the present work.  Despite the difference there are similarities to the present work,
such as the development of velocity integrals similar to the spatially periodic contour dynamics equation we develop for the
periodic patch/layer problem.  Further background on vortex layers may be found in \cite{gargano}.

While we are unaware of other works on the global regularity of unbounded vortex patches for the two-dimensional Euler equations,
the situation is different for the quasi-geostrophic equation.  Rodrigo developed existence theory for a patch which is 
spatially periodic and vertically unbounded in one direction (similarly to a half-space) \cite{rodrigo1}, \cite{rodrigo2}.
More recently Hunter, Shu, and Zhang have studied the related front solutions of the surface quasi-geostrophic equation 
\cite{hunter1}, \cite{hunter2}, \cite{hunter3}.

\subsection{Organization of this paper} We will find many of our calculations much more convenient to perform in the complex plane, yet our results are all real-valued. We describe how to translate back and forth between these settings, largely a matter of notation, in \cref{S:RealComplex}. In \cref{S:PeriodicSolutions} we describe the process of symmetrizing in pairs that is behind the Type 1 solutions, which we explore in \cref{S:Type1}. In \cref{S:Type2} we describe the results of \cite{AfendikovMielke2005,GallaySlijepcevic2014,GallaySlijepcevic2015} that yield Type 2 solutions, and we use those results in \cref{S:Type3} to obtain Type 3 solutions. We show the equivalence of the three types of solution in \cref{S:ThreeTypes}. In \cref{S:VelocityGradient} we give expressions for the velocity gradient in terms of the vorticity, deferring the proofs to \cref{A:VelocityGradient}. We then specialize to vortex patch solutions for Type 1, 2, and 3 solutions, obtaining their contour dynamics equation in \cref{S:CDEPeriodic}, and establishing the global-in-time propogation of the regularity of a vortex patch boundary in \cref{S:BoundaryRegularity}.

%
%
\section{Preliminaries: $\R^2$ and $\C$}\label{S:RealComplex}

\subsection{Real to complex translation}

\noindent Some of our calculations will be more easily performed using complex analysis, though the end results are all real-valued functions. For this we need a means, and a corresponding notation, to switch back and forth between viewing points in the plane as vectors or points in $\C^2$. For this purpose, we will use bold-face letters, such as $\x$ or $\uu$, for quantities that are intrinsically elements of $\R^2$ or vector-valued. We define maps,
\begin{align*}
	\left\{
	\begin{array}{ll}
		\bigvec{} \, \colon \C \to \R^2, \\[4pt]
		\bigvec{x + i y} = (x, y)
	\end{array}
	\right\}
	\text{ and }
	\left\{
	\begin{array}{ll}
		\vecinv{} \colon \R^2 \to \C, \\[4pt]
		\vecinv{(x, y)} = x + i y
	\end{array}
	\right\}.
\end{align*}

For a vector $\x = (x, y)$, we define
\begin{align*}
	\x^\perp := (-y, x).
\end{align*}
Hence, $\x^\perp$ is $\x$ rotated 90 degrees counterclockwise.

\begin{lemma}\label{L:vec}
Let $z, w \in \C$ and $\cdot$ be the usual dot (inner) product of Euclidean vectors. Then
\begin{align}\label{e:vecIds}
	\begin{split}
	\RE \, (zw) &= \vec{\ol{z}} \cdot \vec{w}, \\
	\IM \, (zw) &= -\vec{\ol{z}} \cdot \vec{w}^\perp.
	\end{split}
\end{align}
If $a \in \R$, $z \in \C$,
\begin{align}\label{e:vecProp}
	\bigvec{a z} = a \vec{z}, \quad
	\bigvec{i z} = \vec{z}^\perp, \quad
	\vecinv{\vv^\perp} = i \vecinv{\vv}.
\end{align}
Also, $f$ is analytic in some domain $U$ if and only if $\dv \vec{\ol{f}} = \curl \vec{\ol{f}} = 0$ in $U$, where for any vector field $\vv$,
\begin{align*}
	\dv \vv &:= \pdx{v^1}{x_1} + \pdx{v^2}{x_2}, \quad
	\curl \vv := \pdx{v^2}{x_1} - \pdx{v^1}{x_2}
\end{align*}
are the divergence and (scalar) curl of $\vv$.
\end{lemma}

The boundary integrals we encounter will be real path integrals, but we will sometimes find it useful to transform them to complex contour integrals as in the following lemma:
\begin{lemma}\label{L:ComplexToRealContourIntegrals}
	Let $\bgamma \colon [a, b] \to \C$ be a Lipschitz-continuous path on which the complex-valued
	function $f$ is continuous. Let $\BoldTau$ be the unit tangent vector
	in the direction of $\bgamma$ and $\n$ the associated unit normal,
	with $(\n, \BoldTau)$ in the standard orientation of $(\e_1, \e_2)$.
	Let $C = \image \bgamma$. Then
	\begin{align*}
		\complexint_\bgamma f
			&= \int_C \vec{\ol{f}} \cdot \BoldTau
				+ i \int_C \vec{\ol{f}} \cdot \n.
	\end{align*}
	Here, $\complexint$ is a complex contour integral.
\end{lemma}

Using \cref{L:vec}, it is not hard to rewrite the classical Biot-Savart law in the following hybrid real-complex form:
\begin{theorem}\label{T:BSLawNS2}
	Assume that
	$\omega \in L^1 \cap L^\iny(\R^2)$.
	With $K$ as in \cref{e:BSR2},
	\begin{align}\label{e:BSLawNS2}
		\uu(\x)
			&:= K * \omega(\x)
			= -\bigvec{\frac{i}{2 \pi}
				\int_{\R^2} 
					\frac{\omega(\y)}{\text{ } \ol{\vecinv{\y - \x}} \text{ }} \, d \y}
	\end{align}
	is divergence-free with $\curl \uu = \omega$, and $\uu$ is the unique such velocity
	field in $L^\iny \cap H^1(\R^2)$.
\end{theorem}

\subsection{The cotangent}

\begin{lemma}\label{L:ML}
	For any $z \in \C$ that is not an integer,
	\begin{align*}
		\pi \cot \pi z
			&= \frac{1}{z} + 2 \sum_{n = 1}^\iny \frac{z}{z^2 - n^2}
			= \lim_{N \to \iny} \sum_{n = -N}^N \frac{1}{z + n}.
	\end{align*}
\end{lemma}
\begin{proof}
	For the first equality see, for instance, Equation (11) in Section 5.2.1 of \cite{Ahlfors}.
	The second equality then follows from
	\begin{align*}
		\frac{z}{z^2 - n^2}
			&= \frac{z}{(z - n) (z + n)}
			= \frac{1}{2} \brac{\frac{1}{z - n} + \frac{1}{z + n}}
	\end{align*}
	and summing in pairs, $n$ with $-n$.
\end{proof}

\begin{lemma}\label{L:SumBringsIncot}
	For any $\x, \y \in \R^2$,
	\begin{align*}
		\lim_{N \to \iny} \sum_{n = -N}^N
				\frac{\x + (n, 0)}{\abs{\x + (n, 0)}^2} \cdot \y
			= \pi \bigvec{\cot (\pi \ol{\vecinv{\x}})} \cdot \y.
	\end{align*}
\end{lemma}
\begin{proof}
	Letting $z = \vecinv{\x}$, $w = \vecinv{\y}$, and using \cref{{e:vecIds}}$_1$, we have
	\begin{align*}
		\frac{\x + (n, 0)}{\abs{\x + (n, 0)}^2} \cdot \y
			&= \frac{\RE ((\ol{z} + n) w)}{\abs{z + n}^2}
			= \RE \frac{(\ol{z} + n) w}{\abs{z + n}^2}
			= \RE \frac{w}{z + n}
	\end{align*}
	so
	\begin{align*}
		\lim_{N \to \iny} &\sum_{n = -N}^N
				\frac{\x + (n, 0)}{\abs{\x + (n, 0)}^2} \cdot \y
			= \RE \brac{w \lim_{N \to \iny} \sum_{n = -N}^N \frac{1}{z + n}}
			= \pi \RE (w \cot \pi z) \\
			&= \pi \bigvec{\ol{\cot (\pi \vecinv{\x})}} \cdot \y
			= \pi \bigvec{\cot (\pi \ol{\vecinv{\x}})} \cdot \y,
	\end{align*}
	where we again used \cref{{e:vecIds}}$_1$.
\end{proof}

\subsection{Useful identities}

The identities in \cref{e:sinIdentity,e:TrigIdentities} are easily verifiable; \cref{e:AbramowitzStegun1964} is 4.3.58 of \cite{AbramowitzStegun1964}.

\begin{align}
	\abs{\sin z}^2
		&= \sin^2 x + \sinh^2 y,
			\label{e:sinIdentity} \\
	\cosh 2 x &= 2 \sinh^2 x + 1, \quad
	\cos 2 x = 1 - 2 \sin^2 x,
			\label{e:TrigIdentities} \\
	\cot z
		&= \frac{\sin 2 x - i \sinh 2 y}{\cosh 2y - \cos 2x}.
			\label{e:AbramowitzStegun1964}
\end{align}

\subsection{Lifting paths and domains}\label{S:Lifting}

We will find the need, in the proof of \cref{T:CDEType2}, to apply \cref{L:ComplexToRealContourIntegrals} while integrating in $\Pi$ and apply Cauchy's residue theorem.  This could be done directly by introducing a version of the residue theorem for $\Pi$, which is a (flat) analytic manifold. Alternately, we can transform integrals in $\Pi$ to integrals of $x_1$-periodic functions in $\C$ by \textit{lifting} the domain $\Omega$ in $\Pi$ to a suitable domain $\widetilde{\Omega}$ in $\C$. Our main tool for doing this is the lifting of paths from a topological space to a covering space.

Defining
\begin{align*}
	\Cover \colon \C \to \Pi, \quad
	\Cover(x_1 + i x_2) = x_1 - \lfloor x_1 + \tfrac{1}{2} \rfloor + i x_2,
\end{align*}
we see that $(\C, \Cover)$ is a covering space of $\Pi$ (see Section IX.7 of \cite{Co1978}, for instance). This will allow us to lift a path in $\Pi$ to a path in $\R^2$ or $\C$.

\begin{remark}
	$\Rep(f)(\x) = f(p(\x))$, though we do not make direct use of this.
\end{remark}

\begin{definition}\label{D:Lifting}
	A path in the topological space $X$ is a continuous map from an interval $I$ to $X$.
	The path $\widetilde{\bgamma}$ in $\C$ is a \textit{lift} or \textit{lifting} of the path
	$\bgamma$ in $\Pi$ if $\Cover \circ \widetilde{\bgamma} = \bgamma$.
\end{definition}

\begin{lemma}\label{L:Lifting}
	Let $\bgamma$ be a finite length continuous path in $\Pi$ with initial point $\x_0$.
	For any $\widetilde{\x}_0 \in \Cover^{-1}(\x_0)$, there exists a unique lifting
	$\widetilde{\bgamma}$ with initial point $\widetilde{\x}_0$.
\end{lemma}
\begin{proof}
	This is a classical result; see, for instance, Corollary IX.7.5 of \cite{Co1978}.
\end{proof}

This lifting allows us to relate path integrals in $\Pi$ to lifted path integrals in $\R^2$ or $\C$:

\begin{lemma}\label{L:LiftedIntegral}
	Let $\bgamma$ be a Lipschitz-continuous path in $\Pi$ with a lift $\widetilde{\bgamma}$
	as given by \cref{L:Lifting}. For any any continuous function $f$ on $\Pi$,
	\begin{align*}
		\int_\bgamma f
			= \int_{\widetilde{\bgamma}} f \circ \Cover.
	\end{align*}
	Moreover, the normal vector field $\n$ on $\bgamma$ lifts to itself as does $\BoldTau$;
	that is, $\n(\bgamma(\al)) = \n(\widetilde{\bgamma}(\al))$
	for all $\al$ in the domain of $\bgamma$ (which is the same as the domain
	of $\widetilde{\bgamma}$).
\end{lemma}
\begin{proof}
	Suppose that $\bgamma \colon [a, b] \to \Pi$, in which case also
	$\widetilde{\bgamma} \colon [a, b] \to \C$ with $\Cover \circ \widetilde{\bgamma} = \bgamma$.
	Then
	\begin{align*}
		\int_{\widetilde{\bgamma}} f \circ \Cover
			&= \int_a^b f \circ \Cover(\widetilde{\bgamma}(\al))
				\widetilde{\bgamma}'(\al) \, d \al
			= \int_a^b f (\bgamma(\al))
				\bgamma'(\al) \, d \al
			= \int_\bgamma f.
	\end{align*}
	We used that $\widetilde{\bgamma}'(\al) = \bgamma'(\al)$, since locally
	$\widetilde{\bgamma}$ and $\bgamma$ differ by a constant
	(if we view $\bgamma$ as giving values in $\PiR$).
	This also gives that $\n$ and $\BoldTau$ lift to themselves.
\end{proof}

\cref{L:LiftedIntegral} is not, however, the entire story when we lift the entire boundary of a domain in $\Pi$. An immediate difficulty stems from the ambient space $\Pi$, which is topologically a cylinder, having nontrivial fundamental (and first homology) group $\Z$. Let us say that a closed curve on $\Pi$ wraps around the cylinder $n$ times if it crosses $\set{x_1 = 0}$ (any vertical slice would do) $n$ times counted with sign, positive in one direction, negative in the other (arbitrarily fixing which direction is positive).

A closed path that wraps zero times around the cylinder is homotopic to a point and lifts to a closed path in $\C$. A closed path that wraps around the cylinder $n$ times, however, will lift by \cref{L:Lifting} to a non-closed path in $\C$ that contains $\abs{n} + 1$ points of $\x_0 + \Cal{L}$, where we define here and for future use,
\begin{align}\label{e:L}
	\Cal{L}
		:= \set{\Z} \times \set{0}, \quad
	\Cal{L}^*
		:= \Cal{L} \setminus (0, 0),
\end{align}
treated as subsets of $\R^2$ or of $\C$. Since we are lifting paths that are boundary components, they will always be closed in $\Pi$, but can wrap only $0$ or $\pm 1$ times around the cylinder else they would of necessity self-intersect. 

\begin{figure}[ht]

\IncludeFigure{0.14}{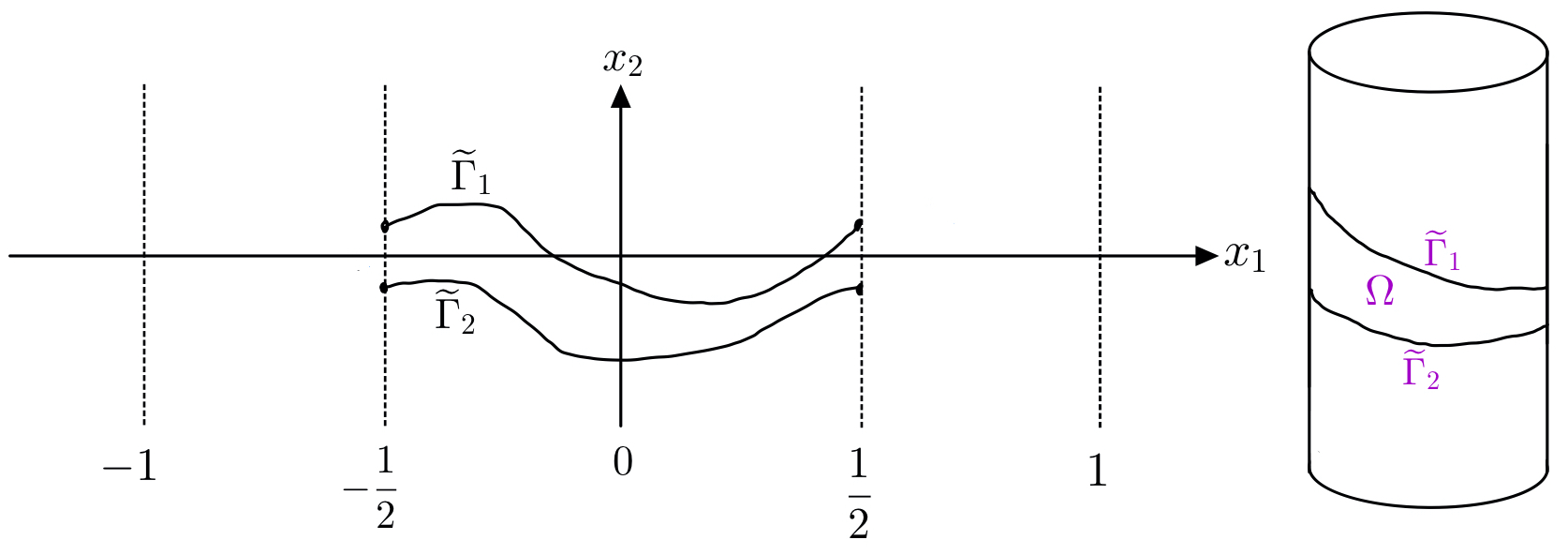}
		
\vspace{-1.5em}
	
\caption{\textit{Lifting of $\prt \Omega$ with base points at
$x_1 = - \tfrac{1}{2}$}}\label{f:BoundaryLift}
\end{figure}

\Cref{f:BoundaryLift} shows an example of a domain $\Omega$ in $\Pi$ having two boundary components $\Gamma_1$, $\Gamma_2$ which lift to non-closed paths $\widetilde{\Gamma}_1$, $\widetilde{\Gamma}_2$. To make a domain from these paths, we could connect $\widetilde{\Gamma}_1$, $\widetilde{\Gamma}_2$ with vertical paths at $x_1 = - \tfrac{1}{2}$ and $x_1 = \tfrac{1}{2}$, oppositely oriented, so that the four paths together form the boundary of a lifted domain $\widetilde{\Omega}$.

Equivalently, and in a manner more easily generalizable, we cut the cylinder $\Pi$ vertically\footnote{In pathological cases, we would have to perturb this cut to avoid producing an infinite number of boundary components, but we will not explore this issue further.} along the line $\ell = \set{x_1 = \pm \tfrac{1}{2}}$, which in effect means we view $\Pi$ in the form suggested in \cref{e:Pi}. For any line segment formed by $\ell \cap \prt \Omega$ we introduce oppositely oriented paths; together, the lifted components of $\prt \Omega$ and these paths, properly oriented, form the boundary components of the lifted domain $\widetilde{\Omega}$.

In lifting these components and paths, however, we need to insure compatible initial points for the paths. To do this, fix any $\x_0$ in $\Omega$. Let $\y$ be any point in $\Omega$ and let $\bgamma_\y$ be a path connecting $\x_0$ to $\y$. Being a domain, $\Omega$ is path-connected so this is always possible. By \cref{L:Lifting}, there is a unique lifting $\widetilde{\bgamma}_\y$ of $\bgamma_\y$ with initial point $\widetilde{\x}_0$. Then $\widetilde{\Omega} := \cup_{\y \in \Omega} \, \widetilde{\bgamma}_\y$ is the desired lifting of $\Omega$.

Lifted in this way, we have the following lemma:

\begin{lemma}\label{L:LiftedBoundaryIntegral}
	Let $\Omega$ be a bounded domain in $\Pi$ and let $\widetilde{\Omega}$ be the
	lifted domain as described above.
	Let $\bgamma$ be a parameterization of $\prt \Omega$ and $\widetilde{\bgamma}$
	a parameterization of $\prt \widetilde{\Omega}$.
	Let $f$ be any continuous complex-valued function.
	Then
	\begin{align*}
		\complexint_\bgamma f
			= \complexint_{\widetilde{\bgamma}} f \circ \Cover, \quad
		\int_{\prt \Omega} \vec{\ol{f}} \cdot \BoldTau
			= \int_{\prt \widetilde{\Omega}} (\vec{\ol{f}} \circ \Cover) \cdot \BoldTau, \quad
		\int_{\prt \Omega} \vec{\ol{f}} \cdot \n
			= \int_{\prt \widetilde{\Omega}} (\vec{\ol{f}} \circ \Cover) \cdot \n.
	\end{align*}
\end{lemma}
\begin{proof}
	Follows from \cref{L:LiftedIntegral}, since the cuts introduce integrals that
	cancel in pairs.
\end{proof}

%
%
\section{Periodized functions and Biot-Savart kernels}\label{S:PeriodicSolutions}

\noindent
\begin{definition}\label{D:SymVel}
	Let $\omega \in L^1 \cap L^\iny(\R^2)$.
	We say that the velocity field $\uu$ is obtained by symmetrizing in pairs (about $0$)
	if, letting $\omega^{(n)}(\x) = \omega(\x + (n, 0))$, we have
	\begin{align*}
		\uu
			= \KSym[\omega]
			:= K * \omega
				+ \sum_{n = 1}^\iny K*\pr{\omega^{(-n)} + \omega^{(n)}}.
	\end{align*}
\end{definition}

\begin{definition}\label{D:S}
	Let $S = S(\R^2)$ be the Serfati space of bounded, divergence-free vector fields on
	$\R^2$ having bounded vorticity with norm,
	\begin{align*}
		\norm{\uu}_S
			:= \norm{\uu}_{L^\iny(\R^2)} + \norm{\curl \uu}_{L^\iny(\R^2)}.
	\end{align*}
	We define $S(\Pi)$ similarly.
\end{definition}

\begin{remark}
	As shown in (2.11) of \cite{GallaySlijepcevic2015}, for any $\omega \in L^\iny(\Pi)$
	there is a divergence-free vector field $\uu$ in $L^\iny(\Pi)$ and so
	in $S(\Pi)$ for which $\curl \uu = \omega$. $S(\R^2)$ is very different,
	for there is no known general condition on $\omega \in L^\iny(\R^2)$ alone that guarantees
	a $\uu$ in $L^\iny(\R^2)$.
\end{remark}

\begin{prop}\label{P:SymConv}
	For $\omega \in L^\iny_c(\R^2)$,
	let $\uu = \KSym[\omega]$ as in \cref{D:SymVel}.
	Then $\uu \in S(\R^2)$ with
	$\curl \uu = \curl \KSym[\omega] = \Rep (\omega)$.
	Further,
	\begin{align}\label{e:KInf}
		\uu
			&= \KSym[\omega]
			= K_\iny * \omega, \quad
		K_\iny(\x)
			:= \bigvec{-\frac{i}{2} \cot \pi \ol{\vecinv{\x}}}
			= \frac{1}{2} \brac{\bigvec{\cot \pi \ol{\vecinv{\x}}}}^\perp,
	\end{align}
	where we note that $K_\iny$ is periodic in $x_1$ with period $1$
	as is $\uu$. We also have,
	\begin{align}\label{e:KH}
		K_\iny(\x)
			&= K(\x) + H(\x),
	\end{align}
	where $H$ is harmonic on $\R^2 \setminus \Cal{L}^*$, where
	$\Cal{L}^*$ is defined in \cref{e:L}.
\end{prop}
\begin{proof}
	Applying \cref{T:BSLawNS2}, we have
	\begin{align*}
		I_n
			&:= K*\pr{\omega^{(-n)} + \omega^{(n)}}(\x)
			= -\bigvec{\frac{i}{2 \pi}
				\int_{\R^2} 
					\brac{
						\frac{\omega(\y)}{\ol{\vecinv{\y - \x}} - n}
					}
					+
					\brac{
						\frac{\omega(\y)}{\ol{\vecinv{\y - \x}} + n}
					}
					\, d \y
					},
	\end{align*}
	so
	\begin{align*}
		\vecinv{I_n}
			&= -\frac{i}{2 \pi} \int_{\R^2}
					2 \frac{\ol{\vecinv{\y - \x}}}
						{(\ol{\vecinv{\y - \x}})^2 - n^2} \omega(\y) \, d \y.
	\end{align*}
	From \cref{D:SymVel} with \cref{L:ML}, then (the compact support of $\omega$ allows
	us to interchange integration and summation),
	\begin{align*}
		\vecinv{\uu(\x)}
			&= -\frac{i}{2} \int_{\R^2}
					\cot (\pi (\ol{\vecinv{\y - \x}}))
					\omega(\y) \, d \y,
	\end{align*}
	and \cref{{e:KInf}} follows from \cref{e:vecProp}.
	Since the singularity of $\cot (\pi z)$ at $z = 0$ is like $1/(\pi z)$
	and $\omega$ is compactly supported, we see that
	the above integral lies in $L^\iny(\R^2)$. Since the curl of each $I_n$ is
	$\omega^{(-n)} + \omega^{(n)}$ while its divergence is zero
	and the sum converges absolutely and uniformly,
	we know that $\dv \uu = 0$ and $\curl \uu = \Rep (\omega)$.
	
	But 
	$
		\cot z
			= \frac{1}{z} + h(z)
	$
	on $\C \setminus \Cal{L}^*$,
	where $h$ is analytic. From this \cref{e:KH} follows. 
\end{proof}

\begin{prop}
	If $\omega_1 \sim \omega_2$ in $L^\iny_c(\R^2)$ as in \cref{D:Equivalence}
	then $K_\iny * \omega_1 = K_\iny * \omega_2$.
\end{prop}
\begin{proof}
	For any $\omega \in L^\iny_c(\R^2)$,
	\begin{align*}
		K_\iny * &\Per(\omega)(\x)
			= \int_\Pi K_\iny(\x - \y) \Per(\omega)(\y) \, d\y
			= \int_\PiR K_\iny(\x - \y) \sum_{n \in \Z} \omega(\y - (n, 0)) \, d\y \\
			&= \int_\PiR \sum_{n \in \Z} K_\iny(\x - (\y - (n, 0))) \omega(\y - (n, 0)) \, d\y \\
			&= \sum_{n \in \Z} \int_\PiR K_\iny(\x - (\y - (n, 0))) \omega(\y - (n, 0)) \, d\y
			= \sum_{n \in \Z} \int_{\PiR - (n, 0)} K_\iny(\x - \y) \omega(\y) \, d\y \\
			&= \int_{\R^2} K_\iny(\x - \y) \omega(\y) \, d\y
			= K_\iny * \omega(\x).
	\end{align*}
	We were able to interchange the integral and sum here because for any fixed $\x$,
	the compact support of $\omega$ makes all but a finite number of terms in the
	sum zero.
	Hence, if $\omega_1 \sim \omega_2$ then
	$K_\iny * \omega_1 = K_\iny * \Per(\omega) = K_\iny * \omega_2$.
\end{proof}

We will see in \cref{S:BSStrip} that $K_\iny$ also serves as the Biot-Savart kernel on $\Pi$.

%
%
\section{Type 1: Periodized solutions}\label{S:Type1}

\noindent We review here results, obtained variously in \cite{Serfati1995Bounded,AKLL2015,TaniuchiEtAl2010,Kelliher2015}, on bounded vorticity, bounded velocity solutions to the 2D Euler equations in $\R^2$.

Let $a_R$ be a radial cutoff function: $a_R(\cdot) = a(\cdot/R)$ for any $ R > 0$, where $a \in C^\iny(\R^2)$ is radially symmetric and equal to $1$ in a neighborhood of the origin. For definitiveness, we will assume that $a \equiv 1$ on $B_1(0)$, $a \equiv 0$ on $B_2(0)^C$, and $\abs{a} \le 1$ on $\R^2$.

\begin{theorem}[\cite{Kelliher2015}]\label{T:Uinf}
Any weak solution to the Euler equations (Eulerian or Lagrangian) with $\uu \in L^\iny(0, T; S) \cap C([0, T] \times \R^2)$ having vorticity $\omega$ with $\uu(0) = \uu^0$, $\omega(0) = \omega^0$, must satisfy, for some $\UU_\iny \in C([0, T])^2$, the Serfati identity,
	\begin{align}\label{e:SerfatiId}
		\begin{split}
			u^j(t&) - (u^0)^j
				= U^j_\iny(t)
					+ (a K^j) *(\omega(t) - \omega^0) \\
				&\qquad- \int_0^t \pr{\grad \grad^\perp \brac{(1 - a) K^j}}
				\stardot (u \otimes u)(s)
							\, ds,
		\end{split}
	\end{align}
$j = 1, 2$, and the renormalized Biot-Savart law,
\begin{align}\label{e:BSRenorm}
	\uu(t) - \uu^0
		&= \UU_\iny(t) + \lim_{R \to \iny} (a_R K) * (\omega(t) - \omega^0)
\end{align}
on $[0, T] \times \R^2$. Furthermore, the corresponding pressure is of the form,
\begin{align}\label{e:pqPressure}
	p(t, \x)
		= -\UU_\iny'(t) \cdot \x + q(t, \x),
\end{align}
where $q$ grows sublinearly at infinity.
\end{theorem}

\cref{T:Uinf} characterizes solutions to the 2D Euler equations that have bounded vorticity and bounded velocity: their existence and uniqueness under the condition that \cref{e:SerfatiId} holds is shown, for $\UU_\iny \equiv 0$, in \cite{Serfati1995Bounded} and elaborated on in \cite{AKLL2015}, their extension to a general $\UU_\iny$ being a simple matter. Uniqueness under the assumption of sublinear growth of the pressure is established in \cite{TaniuchiEtAl2010}.

Combining these results leads to the following:
\begin{theorem}\label{T:Type1}
	Let $\uu^0 \in S(\R^2)$ and set $\omega^0 = \curl \uu^0$.
	There exists a solution $(\uu, p)$ to the 2D Euler equations
	with $\uu \in L^\iny(0, T; S) \cap C([0, T] \times \R^2)$
	having initial velocity $\uu^0$. Existence and uniqueness hold
	if we require that the solution satisfy any one (and hence all) of
	\crefrange{e:SerfatiId}{e:pqPressure} with $\UU_\iny \equiv 0$.
\end{theorem}

%
%
\section{Type 2: Solutions in an infinite periodic strip}\label{S:Type2}

\noindent
Let $\BUC(\Pi)$ be the space of bounded, uniformly continuous functions, noting that any vector field in $S(\Pi)$ lies in $\BUC(\Pi)$. Well-posedness of solutions to the Navier-Stokes equations for initial velocity in $\BUC(\Pi)$ was established by Afendikov and Mielke in \cite{AfendikovMielke2005}. Building on this, Gallay and Slijep\v{c}evi\'{c} in \cite{GallaySlijepcevic2015} (and see the comments in \cite{Gallay2017}) obtained improved bounds for the case where the initial velocity lies in $S(\Pi)$, having established properties of the pressure in \cite{GallaySlijepcevic2014}. These works are for the Navier-Stokes equations, but as the authors point out, the pertinent estimates are uniform in small viscosity and hold for solutions to the Euler equations as well (by repeating the argument with the viscous terms missing or by using known vanishing viscosity results).

In \cref{T:Type2} we give the well-posedness result as derived from \cite{AfendikovMielke2005,GallaySlijepcevic2014,GallaySlijepcevic2015}, but for this we need to first explore some aspects of the analysis in these references.

\subsection{Biot-Savart kernels}\label{S:BSStrip}

The authors of \cite{AfendikovMielke2005,GallaySlijepcevic2014,GallaySlijepcevic2015} orient their periodic strip (infinite cylinder) horizontally and $S^1$ is, in effect, parametrized from $0$ to $1$. Let $(x_1', x_2')$ be the coordinates for the horizontal strip of \cite{AfendikovMielke2005,GallaySlijepcevic2014,GallaySlijepcevic2015}, while we will keep $(x_1, x_2)$ for our vertical strip. Rotating the horizontal strip 90 degrees counterclockwise induces the change of variables,
\begin{align*}
	x_1' \mapsto x_2, \quad x_2' \mapsto -x_1.
\end{align*}

The Biot-Savart kernel on $\Pi$ used in \cite{AfendikovMielke2005} and (2.7) of \cite{GallaySlijepcevic2014} is $\grad^\perp G$, where
\begin{align*}
	G(x_1', x_2')
		:= \frac{1}{4 \pi}
			\log \pr{2 \cosh(2 \pi x_1') - 2 \cos(2 \pi x_2')}
\end{align*}
is the Green's function for the Dirichlet Laplacian on $\Pi$. In $(x_1, x_2)$ variables,
\begin{align}\label{e:Green}
	G(x_1, x_2)
		:= \frac{1}{4 \pi}
			\log \pr{2 \cosh(2 \pi x_2) - 2 \cos(2 \pi x_1)}.
\end{align}

\begin{lemma}\label{L:KInfgradperG}
	We have $K_\iny = \grad^\perp G$. Moreover,
	$G(\x)
			= (2 \pi)^{-1}
				\log \KInfFunc(\x)
	$,
	where
	\begin{align}\label{e:rho}
		\KInfFunc(\x)
			&:= \pr{\sin^2(\pi x_1) + \sinh^2(\pi x_2)}^{\frac{1}{2}}.
	\end{align}
\end{lemma}
\begin{proof}
From \cref{e:TrigIdentities},
$2 \cosh(2 \pi x_2) - 2 \cos(2 \pi x_1) = 4 \rho(\x)^2$, gives our alternate expression for $G$
(noting that the Green's function on $\Pi$ is unique up to an additive constant).
From \cref{e:Green,e:AbramowitzStegun1964}, we have
\begin{align}\label{e:KInfRealForm}
	\grad^\perp G(x_1, x_2)
		= \frac{1}{2 \pi}
			\frac{(- \pi \sinh (2 \pi x_2), \pi \sin(2 \pi x_1))}
			{\cosh(2 \pi x_2) - \cos(2 \pi x_1)}
		= \frac{1}{2}
			\bigvec{\cot(\pi \ol{z})}^\perp,
\end{align}
matching the expression for $K_\iny$ in \cref{e:KInf}. Here, we used \cref{e:AbramowitzStegun1964}.
\end{proof}

\begin{lemma}\label{L:loglogDiff}
	The function $\log \KInfFunc(\x) - \log \abs{\x}$ is harmonic on
	$\R^2 \setminus \Cal{L}^*$,
	where $\KInfFunc$ is defined in \cref{e:rho}.
\end{lemma}
\begin{proof}
	Letting $z = \vecinv{\x}$, we have, using \cref{e:sinIdentity},
	\begin{align*}
		\log \KInfFunc&(\x) - \log \abs{\x}
			= \frac{1}{2} \log \abs[\bigg]{\frac{\KInfFunc(\x)^2}{\abs{\x}^2}}
			= \frac{1}{2} \log \abs[\bigg]{\frac{\sin z}{z}}^2
			= \log \abs[\bigg]{\frac{\sin z}{z}}
			= \RE \log \frac{\sin z}{z},
	\end{align*}
	which is the real part of a function that is complex analytic
	on $\C \setminus \Cal{L}^*$.	
\end{proof}

\subsection{Mean horizontal values}\label{S:MeanHorizontal}
As observed below Lemma 2.2 of \cite{AfendikovMielke2005}, although $K_\iny \in L^1_{loc}(\Pi)$, $K_\iny^2 \in L^1(\Pi)$ (accounting for the different orientation of the strip). Moreover, convolution with $K_\iny^1$ can be handled by subtracting from $u^2$ its mean horizontal value to give it mean value zero. We summarize here this process as described on page 1748 of \cite{GallaySlijepcevic2014}.

If $\vv(t) \in S(\Pi)$, the mean value of $v^2(t)$ along the horizontal line segment $x_2 = a$ is independent of $a \in \R$, and if $(\vv, p)$ solves the Euler equations on $\Pi$ then it is independent of time as well. Hence, we can define
\begin{align}\label{e:m2}
	m_2(t)
		&= m_2[\vv(t)]
		:= \innp{v^2(t)},
\end{align}
the mean value of $v^2(t)$ along any such horizontal line segment and we will have $\innp{v^2(t)} = \innp{v^2_0}$.

The mean value of $v^1(t)$, however, will depend upon $x_2$, so we write
\begin{align*}
	m_1(t, x_2)
		&= m_1[\vv(t)](x_2)
		:= \int_{-\frac{1}{2}}^{\frac{1}{2}} v^1(t, x_1, x_2) \, d x_1.
\end{align*}

Similarly, we define
\begin{align*}
	\innp{\omega}(t, x_2)
		:= \int_{-\frac{1}{2}}^{\frac{1}{2}} \omega(t, x_1, x_2) \, d x_1
\end{align*}
and $\wh{\omega}(t, x_1, x_2) := \omega(t, x) - \innp{\omega}(t, x_2)$. Also,
\begin{align}\label{e:omegaMeanprt2m1}
	\innp{\omega}(t, x_2)
		&= \innp{\prt_1 u^2 - \prt_2 u^1}(t, x_2)
		= -\innp{\prt_2 u^1}(t, x_2)
		= - \prt_2 m_1(t, x_2).
\end{align}

A form of the Biot-Savart law given in (2.5, 2.6) of \cite{GallaySlijepcevic2014} (suppressing the time variable) is
\begin{align}\label{e:BSGS2014}
	\vv(\x)
		&=
		\begin{pmatrix}
			-m_1(x_2) \\
			m_2
		\end{pmatrix}
		+
		\int_{-\iny}^\iny \int_{-\frac{1}{2}}^{\frac{1}{2}}
			K_\iny(\x - \y) \wh{\omega}(\y) \, d y_1 \, d y_2.
\end{align}
We note here that in transforming from the expression as written in \cite{GallaySlijepcevic2014}, a velocity $(v^1, v^2)$ in $(x_1', x_2')$ becomes $(v^2, - v^1)$ in $(x_1, x_2)$, which accounts for the minus sign in $-m_1(x_2)$.

\subsection{Type 2 solutions}
We can now summarize the known result we need for Type 2 solutions:

\begin{theorem}[\cite{AfendikovMielke2005,GallaySlijepcevic2014,GallaySlijepcevic2015}]\label{T:Type2}
	For $\vv^0 \in S(\Pi)$ with $\innp{v^2_0} = 0$
	there exists a unique solution $(\vv, q)$ to the Euler
	equations,
	\begin{align}\label{e:EType2}
		\begin{cases}
			\prt_t \vv + \vv \cdot \grad \vv + \grad q = 0
				&\text{in } [0, \iny) \times \Pi, \\
			\dv \vv = 0
				&\text{in } [0, \iny) \times \Pi, \\
			\vv(0) = \vv^0	
				&\text{in } \Pi
		\end{cases}
	\end{align}
	for which $m_2(t) \equiv 0$
	with $\vv \in C([0, \iny); BUC(\Pi)) \cap L^\iny([0, \iny); S(\Pi))$
	and pressure $q \in W^{1, \iny}([0, \iny) \times \Pi)$.
	The pressure is given by\footnote{$+ 2 K_\iny^2$ is $-\prt_2 G$ in
	(2.8) of \cite{GallaySlijepcevic2014}: we have made the transformation
	from a horizontal to a vertical strip.}
	\begin{align*}
		q
			= - (u^2)^2 + 2 K_\iny^2 * (\omega u^1).
	\end{align*}
	The solutions are Eulerian in velocity and satisfy the vorticity equation.
	Moreover, $\uu$ can be recovered from $\omega$ by the Biot-Savart law 
	as in \cref{e:BSGS2014}.
\end{theorem}

\subsection{Compactly supported vorticity}

As a prelude to obtaining Type 3 solutions, let us consider the special case of Type 2 solutions that we can obtain when the vorticity is compactly supported in $\Pi$. First, we specialize the Biot-Savart law in \cref{e:BSGS2014} to compactly supported vorticity.

\begin{lemma}\label{L:KInfm2m1}
	Let $\vv \in S(\Pi)$ with $\omega := \curl \vv$ compactly supported in $\Pi$.
	Then $m_1(-\iny) + m_1(\iny) \equiv m_2 \equiv 0$
	if and only if $\vv = K_\iny * \omega$.
\end{lemma}
\begin{proof}
	Since $\innp{\omega} = - \prt_2 m_1$, we have
	\begin{align*}
		I^j
			&:= \brac{\int_{-\iny}^\iny \int_{-\frac{1}{2}}^{\frac{1}{2}}
			K_\iny(\x - \y) \innp{\omega}(\y) \, d y_1 \, d y_2}^j \\
			&= \int_{-\frac{1}{2}}^{\frac{1}{2}}
				\int_{-\iny}^\iny
					K_\iny^j((x - x', y - y')) \prt_2 m_1(y') \, d y' \, d x'.
	\end{align*}
	\cref{L:IntK}, below, gives that $I^2 = 0$.

	We now consider $I^1$. Because $\omega$ is compactly supported within some
	$[-1/2, 1/2] \times [-R_0, R_0]$, so,
	too, are $\innp{\omega}$ and then, by \cref{e:omegaMeanprt2m1}, $\prt_2 m_1$.
	Choose $\varphi \in C_C^\iny(\R)$
	equal to 1 on $[-R, R]$ and equal to zero outside $[-R + 1, R + 1]$ where
	we will choose $R \ge R_0$ more precisely later.
	Let $m_1^\eps = \eta_\eps * m_1$, where $\eta_\eps$ is a (compactly supported)
	Friedrich's mollifier.
	As in \cite{AfendikovMielke2005}, we treat $K_\iny^1$ as a distribution on $\Pi$
	with $\varphi m_1^\eps$ a test function. Since also $K_\iny^1 \in L^1_{loc}(\Pi)$,
	we have, for fixed $\x$,
	\begin{align*}
		I^1
			&= \lim_{\eps \to 0}
				\int_{-\frac{1}{2}}^{\frac{1}{2}}
				\int_{-\iny}^\iny
					K_\iny^1((x - x', y - y')) \varphi(y) \prt_2 m_1^\eps(y') \, d y' \, d x' \\
			&= \lim_{\eps \to 0}
				K_\iny^1 * (\varphi \prt_2 m_1^\eps)
			= \lim_{\eps \to 0}
				K_\iny^1 * \prt_2 (\varphi m_1^\eps)
					- \lim_{\eps \to 0}
				K_\iny^1 * (\prt_2 \varphi \, m_1^\eps).
	\end{align*}
	Now,
	\begin{align*}
		\prt_2 K_\iny^1
			&= - \prt_2^2 G
			= - \Delta G + \prt_1^2 G
			= - \delta + \prt_1^2 G,
	\end{align*}
	where $G$ is the Green's function for the Dirichlet Laplacian on $\Pi$ as in
	\cref{e:Green} and $\delta$ is the Dirac delta function on $\Pi$. Hence,
	\begin{align*}
		\prt_2 (\varphi m_1^\eps)
			&= m_1^\eps(x_2)
				- \int_{-\iny}^\iny \int_{-\frac{1}{2}}^{\frac{1}{2}}
					\prt_1^2 G((x - x', y - y')) \, dx' m_1^\eps(y') \, d y'
			= m_1^\eps(x_2),
	\end{align*}
	where the integral vanishes after integrating by parts, since $G$
	is periodic in $x_1$. Hence,
	\begin{align*}
		I^1
			&= m_1(x_2)
				- \lim_{\eps \to 0} K_\iny^1 * (\prt_2 \varphi m_1^\eps),
	\end{align*}
	and this equality holds regardless of our choice of $R \ge R_0$.
	Therefore, if we can evaluate $K_\iny^1 * (\prt_2 \varphi m_1^\eps)$ in the limit
	as $R \to \iny$, it will be its common value for all $R \ge R_0$.
	
	We see from \cref{e:KInfRealForm} that $K_\iny^1(x - y) \to \pm 1/2$ as $y_2 \to \pm \iny$
	and $\prt_2 K_\iny^1(x - y) \to 0$ as $y_2 \to \pm \iny$,
	so
	\begin{align*}
		\lim_{R \to \iny} &K_\iny^1 * (\prt_2 \varphi m_1^\eps)
			= \lim_{R \to \iny} \pr{\int_{-R - 1}^{-R} + \int_R^{R + 1}}
				\prt_2 \varphi  K_\iny^1(x - y) m_1^\eps \\
			&= \lim_{R \to \iny} \brac{
					(K_1^\iny m_1^\eps)(-R) - (K_1^\iny m_1^\eps)(R)  
				}
				- \lim_{R \to \iny} \pr{\int_{-R - 1}^{-R} + \int_R^{R + 1}}
				\varphi  \prt_2 K_\iny^1(x - y) m_1^\eps \\
			&= - \frac{1}{2} \lim_{R \to \iny}  \brac{m_1^\eps(-R) + m_1^\eps(R)}.
	\end{align*}
	We also used here that
	$\prt_2 m_1^\eps = - \eta_\eps * \prt_2 \innp{\omega} = 0$ for $R \ge R_0$.
	Since this limit gives the value for all $R \ge R_0$, we can take $\eps \to 0$
	to conclude that
	\begin{align*}
		I_1
			&= m_1(x_2) + \frac{1}{2} \brac{m_1(-\iny) + m_1(\iny)}.
	\end{align*}
	
	Returning to \cref{e:BSGS2014}, then, we see that
	\begin{align}\label{e:vIsKInfConvPlusTerm}
		\vv(t, \x)
			&=
			\frac{1}{2} \begin{pmatrix}
				m_1(-\iny) + m_1(\iny) \\
				2 m_2
			\end{pmatrix}
			+ (K_\iny * \omega(t)(\x).
	\end{align}
	This shows that $m_1(-\iny) + m_1(\iny) \equiv 0$ and $m_2 \equiv 0$
	if and only if $\vv = K_\iny * \omega$.
\end{proof}

\begin{cor}\label{C:BSonPi}
	Let $\omega \in L^\iny_c(\Pi)$. Then $\vv = K_\iny * \omega$ is the unique
	element in $S(\Pi)$ for which $\curl \vv = \omega$, 
	$m_2[\vv] = 0$, and
	$m_1[\vv](-\iny) + m_1[\vv](\iny) = 0$.
\end{cor}

\begin{prop}\label{P:KInfBSIsSameAsType2}
	Assume that $\omega^0 \in L^\iny_c(\Pi)$,
	$\vv^0 = K_\iny * \omega^0$, and $\vv$ is a Type 2 solution as in \cref{T:Type2}
	with $\vv$ given by \cref{e:BSGS2014}.
	Then $\vv(t) = K_\iny * \omega(t)$ for all $t$.
\end{prop}
\begin{proof}
	It follows from \cref{L:KInfm2m1} that $m_1(0, -\iny) + m_1(0, \iny) = 0$.
	But as observed following
	(2.11) of \cite{GallaySlijepcevic2014},
	$
		\prt_t m_1 = - \innp{u^2 \omega},
	$
	which we note vanishes for all sufficiently large $x_2$ because of the compact
	support of $\omega$. Hence, $m_1(t, -\iny) + m_1(t, \iny) = 0$ for all $t$.
	We conclude from \cref{e:vIsKInfConvPlusTerm} that $\vv(t) = K_\iny * \omega(t)$
	for all $t$.
\end{proof}

We used \cref{L:IntK} in the proof of \cref{L:KInfm2m1}, above.
\begin{lemma}\label{L:IntK}
	For all $y \in \R$, $K_\iny^1(x_1, x_2)$ is even in $x_1$ and odd in $x_2$, while
	$K_\iny^2(x_1, x_2)$ is odd in $x_1$ and even in $x_2$.
\end{lemma}
\begin{proof}
	This follows directly from \cref{e:KInfRealForm}, since $\grad^\perp G = K_\iny$.
\end{proof}

%
%
\section{Type 3: Solutions with a periodized kernel}\label{S:Type3}

\noindent

\begin{theorem}\label{T:Type3}
	Let $\omega^0 \in L^\iny_c(\R^2)$.
	There exists a solution $\mu$ to
	\begin{align*}
		\begin{cases}
			\prt_t \mu + \ww \cdot \grad \mu = 0
				&\text{in } [0, \iny) \times \R^2, \\
			\ww = K_\iny*\mu
				&\text{in } [0, \iny) \times \R^2, \\
			\mu(0) = \omega^0
				&\text{in } \R^2.
		\end{cases}
	\end{align*}
	Moreover, $\curl \ww = \Per (\mu)$,
	and $\ww \in L^\iny(0, T; S) \cap C([0, T] \times \R^2)$ is the unique solution to
	\begin{align}\label{e:Type3Vel}
		\begin{cases}
			\prt_t \ww + \ww \cdot \grad \ww + \grad r = 0
				&\text{in } [0, \iny) \times \R^2, \\
			\dv \ww = 0
				&\text{in } [0, \iny) \times \R^2, \\
			\ww(0) = K_\iny * \mu^0			
				&\text{in } \R^2,
		\end{cases}
	\end{align}
	with the uniqueness criteria being that $r$ is periodic. Finally,
	$r \in L^\iny([0, T] \times \R^2)$.
\end{theorem}
\begin{proof}
	From \cref{P:SymConv} we know that $K_\iny * \omega^0 \in L^\iny(\R^2)$ and is
	periodic in $x_1$ with period 1; hence, abusing notation, we can set
	$\vv^0 = K_\iny * \omega^0|_\Pi$ and obtain by \cref{T:Type2}
	a unique solution $(\vv, q)$ to \cref{e:EType2}
	for which $q$ is periodic in $x_1$ and $m_2(t) \equiv 0$.
	Since $\curl \vv^0 = \omega^0|_\Pi$ is compactly supported and so
	$\curl \vv$ remains compactly supported for all time,
	we know from \cref{P:KInfBSIsSameAsType2}
	that $\vv = K_\iny*\curl \vv$.
	So letting $\zeta = \curl \vv$,
	we see that
	\begin{align*}
		\begin{cases}
			\prt_t \zeta + \vv \cdot \grad \zeta = 0
				&\text{in } [0, \iny) \times \Pi, \\
			\vv = K_\iny * \zeta
				&\text{in } [0, \iny) \times \Pi, \\
			\zeta(0) = \omega^0				
				&\text{in } \Pi.
		\end{cases}
	\end{align*}
	Setting $\ww = \vv$, $\mu = \zeta$ gives the desired solution of Type 3.
	Moreover, since $q(t)$ is periodic, we can let $r = \Rep(q)$,
	and we obtain a unique solution to \cref{e:Type3Vel}.
\end{proof}

%
%
\section{Three types of solution are equivalent}\label{S:ThreeTypes}

\noindent For certain classes of initial data, our three types of solution are equivalent. The equivalence of Type 1 and Type 2 holds for a broader class, so we first prove it in \cref{T:TwoSolutions}. The equivalence of the third type holds for initial data in $L^\iny_c(\R^2)$, as we show in \cref{T:ThreeSolutions}. This includes vortex patch data, our application in \cref{S:CDEPeriodic}.

\begin{theorem}\label{T:TwoSolutions}
	Let $\vv^0 \in S(\Pi)$ and periodize it
	to give $\uu^0 = \Rep(\vv^0) \in S(\R^2)$. Let $(\uu, p)$ be the solution of Type 1 with initial
	velocity $\uu^0$ given by \cref{T:Type1} and let $(\vv, q)$ the solution of Type 2
	with initial velocity $\vv^0$ given by \cref{T:Type2}. Then
	$\Rep (\vv) = \uu$.
	\end{theorem}
\begin{proof}
	We have $\curl \vv(0) = \curl \uu^0|_\Pi$, where
	we abuse notation somewhat.
	From \cref{T:Type2}, we have  a pressure $q$ with $q(t) \in L^\iny(\Pi)$
	for which
	\begin{align}\label{e:vPeriodicVelSol}
		\begin{cases}
			\prt_t \vv + \vv \cdot \grad \vv + \grad q = 0
				&\text{in } [0, \iny) \times \Pi, \\
			\dv \vv = 0
				&\text{in } [0, \iny) \times \Pi, \\
			\vv(0) = \vv^0
				&\text{in } \Pi.
		\end{cases}
	\end{align}
	
	Since $\Rep(\vv)$ and $\Rep(q)$ are $x_1$-periodic with period 1, we can set
	$\wv = \Rep (\vv)$ and $\wq = \Rep (q)$, and both will lie in $L^\iny([0, T] \times \R^2)$
	with $\curl \wv(t) = \Rep (\curl \vv(t))$. Thus,
	$\wv$ is $\vv$ periodized and $\curl \wv$ is $\curl \vv$ periodized, meaning that
	\cref{e:vPeriodicVelSol} in effect holds on $\Pi_p$
	translated by $(n, 0)$ for any integer $n$, so we see that
	\begin{align}\label{e:uR2VelSol}
		\begin{cases}
			\prt_t \wv + \wv \cdot \grad \wv + \grad \wq = 0
				&\text{in } [0, \iny) \times \R^2, \\
			\dv \wv = 0
				&\text{in } [0, \iny) \times \R^2, \\
			\wv(0) = \uu^0
				&\text{ in } \R^2.
		\end{cases}
	\end{align}
	We see that $(\wv, \wq)$ is a solution to the Euler equations on
	$[0, \iny) \times \R^2$.	
	Manifestly, $\wv$, $\curl \wv$, and $\wq$ each lie
	in $L^\iny([0, \iny) \times \R^2)$, being periodic in $x_1$.
	Hence, $\wv$ is a bounded velocity, bounded vorticity solution
	to the Euler equations on $[0, \iny) \times \R^2$.
	Because the pressure $\wq$ grows sublinearly it is, in fact,
	the (unique) Serfati solution (it satisfies the Serfati identity),
	as follows from \cref{T:Type1}.
	Therefore, $\uu = \vv$.
\end{proof}

\begin{theorem}\label{T:ThreeSolutions}
	For $\omega^0 \in L^\iny_c(\R^2)$,
	let $\uu^0 = \KSym[\omega^0]$ be
	obtained by symmetrizing in pairs as in \cref{D:SymVel},
	and let $\vv^0 = K_\iny * \Per(\omega^0)$.
	Let $(\uu, p)$, $(\vv, q)$ be the Type 1, 2 solutions with
	initial velocity $\uu^0$, $\vv^0$ and let $\ww^0$ be the velocity
	field for the Type 3 solution given by  \cref{T:Type3}.
	Then $\Rep (\vv) = \uu = \ww$.
\end{theorem}
\begin{proof}
	\cref{T:TwoSolutions} gives $\Rep (\vv) = \uu$,
	while  $\Rep (\vv) = \ww$ is inherent in the proof of \cref{T:Type3}.
\end{proof}

%
%
\section{The velocity gradient}\label{S:VelocityGradient}

\noindent 
The following expression for $\grad (K * \omega)$ is classical (see, for instance, Proposition 2.20 of \cite{MB2002}):
\begin{lemma}\label{L:graduClassical}
	Assume that $\omega \in L^\iny(\R^2)$ is compactly supported and let
	$\uu = K * \omega$. Then
	\begin{align*}
		\grad \uu(\x)
			&= \omega(\x)
				\begin{pmatrix}
					0 & -1 \\
					1 & 0
				\end{pmatrix}
				+
				\PV \int_{\R^2} \grad K(\x - \y) \omega(\y) \, d \y,
	\end{align*}
	where we can write,
	\begin{align*}
		\grad K(\x)
			&= \frac{1}{2 \pi} \frac{\sigma(\x)}{\abs{\x}^2}, \quad
		\sigma(\x)
			:= \frac{1}{\abs{\x}^2}
			\begin{pmatrix}
				2 x_1 x_2 & x_2^2 - x_1^2 \\
				x_2^2 - x_1^2 & - 2 x_1 x_2
			\end{pmatrix}.
	\end{align*}
\end{lemma}

The analog for the $K_\iny$ kernel is \cref{L:graduTypeI}.

\begin{lemma}\label{L:graduTypeI}
	Assume that $\omega \in L^\iny(\R^2)$ is compactly supported and let
	$\uu = K_\iny * \omega$. Then
	\begin{align*}
		\grad \uu(\x)
			&= \sum_{n \in \Z} \frac{\omega(\x + (n, 0))}{2}
				\begin{pmatrix}
					0 & -1 \\
					1 & 0
				\end{pmatrix}
				+ 
				\PV \int_{\R^2} \grad K_\iny(\x - \y)
					\omega(\y) \, d \y,
	\end{align*}
	where $\KInfFunc$ is as in \cref{e:rho} and
	where we can write,
\begin{align*}
	\grad K_\iny(\x)
		&= \frac{\pi}{2} \frac{\beta(\x)}{\rho(\x)^2},
\end{align*}
where
\begin{align*}
	\beta(\x)
		&= \frac{1}{2 \rho(\x)^2}
			\begin{pmatrix}
				\sin(2 \pi x_1) \sinh(2 \pi x_2) &
					\cos(2 \pi x_1) \cosh(2 \pi x_2) - 1 \\
				\cos(2 \pi x_1) \cosh(2 \pi x_2) - 1 &
					- \sin(2 \pi x_1) \sinh(2 \pi x_2)
			\end{pmatrix}.
\end{align*}
\end{lemma}

\begin{remark}\label{R:gradKInfSize}
	Like $\sigma$, the matrix $\beta$ is symmetric with trace zero.
	Near the origin, $\KInfFunc(\x)^2 \approx \pi^2 \abs{\x}^2$,
	and we can see that
	$\beta(\x) \approx 4 \pi^2 \abs{\x}^2 / (2 \pi^2 \abs{\x}^2) \approx 2 \approx \sigma(\x)$, and so
	$\grad K_\iny(\x) \approx 2 \pi / (2 \pi^2 \abs{\x}^2) \approx 1/(\pi \abs{\x}^2)
	\approx \grad K(\x)$. Also like $\sigma$, $\beta_{11}$
	and $\beta_{22}$ integrate to zero over circles centered at the origin,
	but unlike $\sigma$, neither $\beta_{12}$ nor $\beta_{21}$ integrate
	to zero.
\end{remark}

We have the following immediate corollary of \cref{L:graduTypeI}:
\begin{cor}\label{C:graduType}
	Let $\vv \in S(\Pi)$ with $\omega = \curl \vv$ compactly supported
	and let  $\uu = K_\iny * \omega$. Then
	\begin{align*}
		\grad \uu(\x)
			&= \frac{\omega(\x)}{2}
				\begin{pmatrix}
					0 & -1 \\
					1 & 0
				\end{pmatrix}
				+
				\PV \int_\Pi \grad K_\iny(\x - \y)
					\omega(\y) \, d \y
	\end{align*}
	and $\grad K_\iny$ can be written as in \cref{L:graduTypeI}.
\end{cor}

%
%
\section{Contour Dynamics Equations}\label{S:CDEPeriodic}

\noindent
First we review the Contour Dynamics Equation (CDE) for a classical vortex patch---the characteristic function of a bounded, simply connected domain evolving under the vorticity equation for the Euler equations on all of $\R^2$---then turn to the CDE for Type 2 solutions.

In what follows we use the Lipschitz space $Lip$ and homogeneous Lipschitz space $lip$. On $U \subseteq \R^d$ for $d \ge 1$, we define their semi-norm and norm,
\begin{align*}
	\norm{f}_{lip(U)}
		:= \sup_{x \ne y \in U} \frac{\abs{f(x) - f(y)}}{\abs{x - y}}, \quad
	\norm{f}_{Lip(U)}
		:= \norm{f}_{L^\iny(U)} + \norm{f}_{lip(U)}.
\end{align*}

\subsection{Classical vortex patches}\label{S:ClassicalPatches}

In the classical setting of a vortex patch in $\R^2$, we have \cref{T:ClassicalPatchBSLaw,T:ClassicalPatchCDE}, as in Proposition 8.6 of \cite{MB2002}
and the derivation of the classical CDE that appears before it. 

\begin{theorem}\label{T:ClassicalPatchBSLaw}
  	Let $\bgamma \colon [0, 2 \pi] \to \R^2$ be a $C^1$
	counterclockwise\footnote{In \cite{MB2002}, the patch boundary is parameterized clockwise,
	but $(\BoldTau, \n)$ is in the standard $(\e_1, \e_2)$ orientation; the two resulting sign
	changes between \cite{MB2002} and us cancel, so there is no sign change
	in our expressions.} parameterization
	of the boundary of a bounded, simply connected domain $\Omega$.
	Then
	\begin{align}\label{e:BSLawPatchR2}
		\uu(\x)
			= -\frac{\omega_0}{2 \pi} \int_0^{2 \pi}
				\log \abs{\x - \bgamma(\al)} \prt_\al \bgamma(\al) \, d \al
	\end{align}
	is the unique divergence-free vector field decaying at infinity for
	which $\curl \uu = \omega_0 \CharFunc_\Omega$.
\end{theorem}

Now let us suppose that $\Omega$ is a simply connected bounded domain in $\R^2$ with a $C^{1, \eps}$ boundary. Let $\uu$ be the unique weak solution to the Euler equations with initial vorticity $\omega^0 := \omega_0 \CharFunc_\Omega$ and let $X$ be the flow map for $\uu$. Then we know that the vorticity $\omega(t) = \omega_0 \CharFunc_{\Omega_t}$, where $\Omega_t = X(t, \Omega)$.

Let $\bgamma(0, \cdot)$ be a $C^1$-regular counterclockwise parameterization of $\Gamma = \prt \Omega$. Define a parameterization of $\prt \Omega_t = X(t, \Gamma)$ by $\bgamma(t, \cdot) := X(t, \bgamma(0, \cdot))$. The log-Lipschitz regularity of $\uu(t)$ induces $C^{c(t)}$-regularity of the flow map $X(t, \cdot)$ with $c(t) \in (0, 1)$ and decreasing with time, as in Lemma 8.2 of \cite{MB2002}. This is insufficient regularity to obtain a $C^1$-parameterization of $\prt \Omega_t$, so let us \textit{suppose} that our (classical) solution has $\uu \in C(0, T; lip)$. Then $\bgamma(t, \cdot)$ is a $C^1$-parameterization of $\prt \Omega_t$.

 Since we assumed $\prt \Omega$ is $C^{1, \eps}$, we could give $\bgamma(0, \cdot)$ $C^{1, \eps}$-regularity, but this does not itself ensure that $\bgamma(t, \cdot)$ is $C^{1, \eps}$: proving that is tantamount to establishing the propagation of regularity of the vortex patch boundary.

 \begin{theorem}\label{T:ClassicalPatchCDE}
 	Let $\uu(t, \x)$ be given by \cref{e:BSLawPatchR2} applied with
	$\bgamma(t, \cdot)$; that is,
	\begin{align*}
		\uu(t, \x)
			:= -\frac{\omega_0}{2 \pi} \int_0^{2 \pi}
				\log \abs{\x - \bgamma(t, \al)} \prt_\al \bgamma(t, \al) \, d \al.
	\end{align*}
	Then $\uu$ is a weak solution to the 2D Euler equations
	on $[0, T] \times \R^2$ with $\uu \in C(0, T; Lip)$
	if and only if
	$\bgamma$ is a $C^1([-T, T]; C([0, 2 \pi])) \cap C([-T, T]; C^1([0, 2 \pi]))$
	solution to the contour dynamics equations (CDE),
	\begin{align}\label{e:CDEClassical}
		\diff{}{t} \bgamma(t, \al)
			= -\frac{\omega_0}{2 \pi} \int_0^{2 \pi}
				\log \abs{\bgamma(t, \al) - \bgamma(t, \al')}
					\prt_{\al'} \bgamma(t, \al') \, d \al.
	\end{align}
\end{theorem}

\cref{T:ClassicalPatchBSLaw,T:ClassicalPatchCDE} were expressed for simply connected domains. As pointed out on page 330 of \cite{MB2002}, the only difference for multiply connected domains is that the integrals in \cref{e:BSLawPatchR2,e:CDEClassical} are summed over each component of the boundary.

\begin{theorem}\label{T:ClassicalMultiplyConnected}
	\cref{T:ClassicalPatchBSLaw,T:ClassicalPatchCDE} hold for bounded, multiply
	connected domains if we evaluate and sum each of the boundary integrals over each
	boundary component.
\end{theorem}

We view \cref{e:CDEClassical} as a form of the Euler equations applying specifically to a vortex patch: it comes directly from \cref{e:BSLawPatchR2}, which we view as a form of the Biot-Savart law that recovers the velocity from the vorticity, as it is encoded by $\bgamma$. We work, now, to obtain replacements for these expressions that apply to periodized vortex patches. This is a matter of deriving the CDE for a solution to the Euler equations and showing, conversely, that any solution to the CDE satisfies the Euler equations.

\subsection{Type 2 solutions}

Turning to Type 2 solutions, we make the following assumptions on $\Omega$:

\begin{assumption}\label{A:OmegaPi}
	Assume that $\Omega \subseteq \Pi$ is bounded with a finite number of boundary
	components, $\Gamma_1, \dots, \Gamma_J$, each $C^{1, \eps}$ regular.
\end{assumption}

With $\Omega$ as in \cref{A:OmegaPi}, we let $\uu$ be the unique Type 2 solution having initial vorticity $\omega^0 := \omega_0 \CharFunc_\Omega$ with $m_2 \equiv m_1(t, -\iny) + m_2(t, \iny) \equiv 0$ given by \cref{T:Type2,P:KInfBSIsSameAsType2}
($m_1$, $m_2$ are defined in \cref{S:MeanHorizontal}).
Set
\begin{align*}
	\Omega_t
		:= X(t, \Omega), \quad
	\Gamma_{t, j}
		:= X(t, \Gamma_j),
\end{align*}
noting that because $X(t, \cdot)$ is a homeomorphism of $\R^2$ onto $\R^2$, $\Gamma_{t, j}$ is the $j^{th}$ of the $J$ components of $\prt \Omega_t$. We then define a parameterization $\bgamma_j$ of $\Gamma_{t, j}$ as we parameterized $\prt \Omega_t$ in \cref{S:ClassicalPatches}, setting $\bgamma_j(t, \cdot) := X(t, \bgamma_j(0, \cdot))$. As in that section, a priori, we do not even know that $\bgamma_j(t)$ has $C^1$ regularity for $t > 0$; proving that it has $C^{1, \eps}$ regularity is the ultimate goal (of \cref{S:BoundaryRegularity}).

We show in \cref{T:Type2PatchBSLaw,T:CDEType2} that the analog of \cref{T:ClassicalMultiplyConnected} holds for Type 2 solutions.

\begin{theorem}\label{T:Type2PatchBSLaw}
	Let $\Omega$ be as in \cref{A:OmegaPi},
	and for each $j$, let $\bgamma_j \colon [0, 2 \pi] \to \R^2$ be a $C^1$
	counterclockwise parameterization of the boundary component $\Gamma_j$.
	With $\rho$ as in \cref{e:rho},
	\begin{align}\label{e:BSlogsin}
		\begin{split}
		\uu(\x)
			&= -\frac{\omega_0}{2 \pi} \sum_{j = 1}^J \int_0^{2 \pi}
				\log \KInfFunc (
					\x - \bgamma_j(\al) 
					)
					\prt_\al \bgamma_j(\al) \, d \al
		\end{split}
	\end{align}
	is the unique divergence-free vector field in $S(\Pi)$
	having curl equal to $\omega_0 \CharFunc_\Omega$
	for which $m_2 = 0$ and $m_1(-\iny) + m_1(\iny) = 0$.
\end{theorem}
\begin{proof}
	By \cref{C:BSonPi}, we know that $\uu = K_\iny * \omega$
	is the unique divergence-free vector field in $S(\Pi)$
	having curl equal to $\omega_0 \CharFunc_\Omega$
	for which $m_2 = 0$ and $m_1(-\iny) + m_1(\iny) = 0$.
	Then we have, using \cref{L:KInfgradperG}
 	and parameterizing $\Gamma_{t, j}$ by arc length from 0 to $\ell_j$,
	setting $\y(s) = \bgamma_j(\al(s))$,
	\begin{align*}
		\uu(\x)
			&= K_\iny * \omega(\x)
			= \grad^\perp G * \omega(\x)
			= \frac{\omega_0}{2 \pi} \int_\Omega
				\grad^\perp \log \KInfFunc(\x - \y) \, d \y \\
			&= -\frac{\omega_0}{2 \pi} \int_\Omega
				\grad^\perp_\y \log \KInfFunc(\x - \y) \, d \y
			= -\frac{\omega_0}{2 \pi} \sum_{j = 1}^J \int_0^{\ell_j}
				\log \KInfFunc(\x - \y(s)) (-n^2, n^1) \, d s \\
			&= -\frac{\omega_0}{2 \pi} \sum_{j = 1}^J \int_0^{\ell_j}
				\log \KInfFunc(\x - \y(s)) \BoldTau(s) \, d s
			= -\frac{\omega_0}{2 \pi} \sum_{j = 1}^J \int_0^{2 \pi}
				\log \KInfFunc (
					\x - \bgamma_j(\al) 
					)
					\prt_\al \bgamma_j(\al) \, d \al.
	\end{align*}
	Here $(n^1, n^2) = \n$ and $(-n^2, n^1) = \BoldTau$ (see \cref{L:ComplexToRealContourIntegrals}),
	and we used that
	\begin{align*}
		\prt_\al \bgamma_j(\al) \, d \al
			&= \frac{\prt_\al \bgamma_j(\al)}{\abs{\prt_\al \bgamma_j(\al)}}
				\abs{\prt_\al \bgamma_j(\al)} \, d \al
			= \BoldTau(s) \, ds.
	\end{align*}
	From this, \cref{e:BSlogsin} follows.
\end{proof}

\begin{theorem}\label{T:CDEType2}
	Let $\uu$ be the Type 2 solution described above
	and assume that
	each $\bgamma_j$ is in
	$C^1([-T, T]; C([0, 2 \pi])) \cap C([-T, T]; C^1([0, 2 \pi]))$.
	Then
	\begin{align}\label{e:BSlogsintime}
		\uu(t, \x)
			&= -\frac{\omega_0}{2 \pi} \sum_{j = 1}^J \int_0^{2 \pi}
				\log \KInfFunc (
					\x - \bgamma_j(t, \al) 
					)
					\prt_\al \bgamma_j(t, \al) \, d \al
	\end{align}	
	and lies in $C(0, T; Lip)$.
	Moreover, each $\bgamma_k$ satisfies the CDE,
	\begin{align}\label{e:CDElogsin}
		\begin{split}
		\diff{}{t} \bgamma_k(t, \al)
			&= -\frac{\omega_0}{2 \pi} \sum_{j = 1}^J \int_0^{2 \pi}
				\log \KInfFunc (
					\bgamma_k(t, \al) - \bgamma_j(t, \al') 
					)
				\prt_\al \bgamma_j(t, \al) \, d \al'.
		\end{split}
	\end{align}
	
	Conversely, if each $\bgamma_k$ in
	$C^1([-T, T]; C([0, 2 \pi]))$ $\cap$ $C([-T, T]; C^1([0, 2 \pi]))$
	satisfies \cref{e:CDElogsin}
	then
	$\uu$ given by \cref{e:BSlogsintime} is a Type 2 solution
	with $\uu \in C(0, T; Lip)$ and $m_2 \equiv m_1(t, -\iny) + m_2(t, \iny) \equiv 0$.
\end{theorem}
\begin{proof}
	The forward direction follows directly from \cref{T:Type2PatchBSLaw}.
	
	For the converse, we parallel the proof of
	Proposition 8.6 of \cite{MB2002}, which consists of two steps: (1) Show that
	$\uu$ given by \cref{e:BSlogsintime} is divergence-free with
	$\curl \uu = \CharFunc \Omega_{0, t}$. (2) Show that $\uu$ solves the 2D Euler equations.
			
	To prove (1), let $\uu$ be given by \cref{e:BSlogsintime}.
	Reparameterizing by arc length as in the proof of \cref{T:Type2PatchBSLaw},
	\begin{align*}
		\uu(t, \x)
			&= -\frac{\omega_0}{2 \pi} \sum_{j = 1}^J \int_0^{\ell_j}
				\log \KInfFunc (\x - \y(s))
					\BoldTau(s) ds
			= -\frac{\omega_0}{2 \pi} \sum_{j = 1}^J \int_{\Gamma_{t, j}}
				\log \KInfFunc (\x - \cdot)
					\BoldTau.
	\end{align*}
	To apply $\dv$ and $\curl$ to this expression, we use that
	for a constant vector field
	$\ww$ and scalar function $g$, $\dv (g \ww) = \grad g \cdot \ww$ and
	$\curl (g \ww) = \grad^\perp g \cdot \ww$. Also, letting
	$\vv = (2 \pi)^{-1} \grad^\perp \log \KInfFunc (\x - \cdot)$
	and $f = \ol{\vecinv{\vv}}$, we see that
	\begin{align*}
		\curl \uu(t, \x)
			&= -\frac{\omega_0}{2 \pi} \sum_{j = 1}^J \int_{\Gamma_{t, j}}
				\grad^\perp \log \KInfFunc (\x - \cdot)
					\cdot \BoldTau
			= -\omega_0 \sum_{j = 1}^J \int_{\Gamma_{t, j}}
				\vec{\ol{f}} \cdot \BoldTau
			= -\omega_0 \int_{\prt \Omega_t}
				\vec{\ol{f}} \cdot \BoldTau, \\
		\dv \uu(t, \x)
			&= -\frac{\omega_0}{2 \pi} \sum_{j = 1}^J \int_{\Gamma_{t, j}}
				\grad \log \KInfFunc (\x - \cdot)
					\cdot \BoldTau
			= \frac{\omega_0}{2 \pi} \sum_{j = 1}^J \int_{\Gamma_{t, j}}
				\grad^\perp \log \KInfFunc (\x - \cdot)
					\cdot \n \\
			&= \omega_0 \sum_{j = 1}^J \int_{\Gamma_{t, j}}
				\vec{\ol{f}} \cdot n
			= \omega_0 \int_{\prt \Omega_t}
				\vec{\ol{f}} \cdot n.
	\end{align*}
	
	Up to this point, we have been integrating over paths in $\Pi$ treated as $\R^2 / \Cal{L}$,
	but we wish to apply \cref{L:ComplexToRealContourIntegrals}, which obliges us to
	work in $\C$. To do this, we lift $\Omega_t$ to $\widetilde{\Omega}_t$ as
	described in \cref{S:Lifting}.
	Applying \cref{L:LiftedBoundaryIntegral,L:ComplexToRealContourIntegrals}
	(writing $f$ in place of $f \circ \Cover$ by viewing $f$ as $x_1$-periodic with period $1$)
	gives for all
	$\x$ not lying on $\prt \widetilde{\Omega}_t$ (a set of measure $0$),
	\begin{align*}
		\omega_0 \, \complexint_{\prt \Omega_t} f
			&= \omega_0 \, \complexint_{\prt \widetilde{\Omega}_t} f
			= \omega_0 \int_{\prt \widetilde{\Omega}_t} \vec{\ol{f}} \cdot \BoldTau
				+ i \omega_0 \int_{\prt \widetilde{\Omega}_t} \vec{\ol{f}} \cdot n
			= -\curl \uu(t, \x) + i \dv \uu(t, \x).
	\end{align*}
	But we see from \cref{L:KInfgradperG} that $\vv = K_\iny(\x - \cdot)$ and that
	\begin{align*}
		f
			&= \frac{1}{2} \ol{\vecinv{\bigvec{\cot(\pi \ol{z})}^\perp}}
			= \frac{1}{2} \ol{\vecinv{\bigvec{i \ol{\cot(\pi z)}}}}
			= \frac{1}{2} \ol{i \ol{\cot(\pi z)}}
			= -\frac{i}{2} \cot(\pi z),
	\end{align*}
	where we used \cref{e:vecProp} and the identity $\ol{i \ol{z}} = - i z$.
	The complex meromorphic function $f$ has simple poles at each point in
	$\x + \Cal{L}$ with residue $(-2 \pi)^{-1} i$.
	By the residue theorem, then,
	summing over all points of $\Cal{L}$ lying inside $\prt \widetilde{\Omega}_t$---that is,
	lying in $\widetilde{\Omega}_t$, 
	\begin{align*}
		\omega_0 \, \complexint_{\prt \widetilde{\Omega}_t} f
			&= \RE \pr{
				2 \pi i \omega_0 \sum_n \Res (f, (n, 0))
				}
			= \omega_0 \RE \pr{
				\frac{2 \pi i}{-2 \pi i}
					\sum_n 1
				}
			= -\omega_0 n.
	\end{align*}
	But $\widetilde{\Omega}_t$ can contain at most one point of $\x + \Cal{L}$ else
	the lift given in \cref{S:Lifting} would map $\x$ to more than one point
	in $\C$ (which would mean it is not a lift).
	We see, then, that
	\begin{align*}
		\curl \uu(t, \x)
			&= -\omega_0 \, \complexint_{\prt \Omega_t} f
			= -\omega_0 \, \complexint_{\prt \widetilde{\Omega}_t} f
		 	= \omega_0 \CharFunc_{\Omega_t}(t, \x)
			= \omega(t, \x).
	\end{align*}
	
	We conclude that for all $t \in [0, T]$, $\dv \uu = 0$ and
	$\curl \uu = \omega = \omega_0 \CharFunc_{\Omega_t}$.
	Directly from \cref{e:BSlogsintime}, we know that $\uu \in L^\iny(\Pi)$ 
	and hence $\uu \in S(\Pi)$. It follows from \cref{T:Type2PatchBSLaw} applied
	with $\gamma_j(t, \cdot)$ in place of $\gamma_j$ for any fixed $t$
	that $m_2[\uu(t)] = 0$ and $m_1[\uu(t)](-\iny) + m_1[\uu](\iny) = 0$.
	
	Using (1), the proof of (2) that $\uu$ solves the 2D Euler equations
	on the time interval $[-T, T]$ proceeds
	just as it does in the proof of Proposition 8.6 on page 334 of \cite{MB2002}.
\end{proof}

\begin{remark}\label{R:Patch1and3}

We can view Type 2 solutions as equivalent to Type 1 or 3 solutions by virtue of \cref{T:ThreeSolutions}. For vortex patches it is most natural to start with an $\Omega \in \Pi$ satisfying \cref{A:OmegaPi} and lift it to $\R^2$ as in \cref{S:Lifting} to give $\Omega_0$. It is also possible to start with a domain in $\R^2$, and use it to obtain via the $\Per$ operator a domain in $\Pi$, but there are no simple general conditions to guarantee that the boundary of the domain in $\Pi$ is regular.
\end{remark}

%
%
\section{Regularity of a vortex patch boundary}\label{S:BoundaryRegularity}

\noindent To prove the propagation of regularity of a vortex patch boundary for our Type 1, 2, or 3 solutions, it will be easiest to work with Type 2 solutions, the result then immediately following for the other two types by \cref{T:ThreeSolutions}. We will prove, in \cref{T:PatchType2}, that for Type 2 solutions, the regularity of the boundary of a periodic vortex patch is maintained for all time, as in the classical case.

\begin{theorem}\label{T:PatchType2}
	Let $\Omega$ be as in \cref{A:OmegaPi} and let $\Omega_t = X(t, \Omega)$ for a Type 2 solution.
	Then $\prt \Omega_t$ is $C^{1, \eps}$ for all time. The analogous result holds for
	Type 1 and 3 solutions.
\end{theorem}
\begin{proof}
We describe only how the proof differs from the now classical proof as presented in Chapter 8 of \cite{MB2002}.
There are two main steps to the proof given in \cite{MB2002}: First, show local-in-time existence of a $C^{1, \eps}$ solution to the CDE (based on \cite{BertozziThesis1991}) then show that the solution extends globally in time (based on \cite{ConstantinBertozzi1993}).

\bigskip\noindent\textbf{Local-in-time $C^{1, \eps}$ solutions}: In brief, the first step is to define the function $F$ on the space $B^{1, \eps}$ of closed $C^{1, \eps}$ paths in $\Pi$ by (we have translated this to Type 2 solutions) by
\begin{align*}
	F(\bgamma(\beta))
		:= \frac{\omega_0}{2 \pi} \int_0^{2 \pi}
				\log \KInfFunc (
					\bgamma(\beta) - \bgamma(\al) 
					)
				\prt_\al \bgamma(\al) \, d \al.
\end{align*}
Here, $F$ is as defined for each boundary component separately, we suppress the sums over each boundary component for notational simplicity. First show that $F \colon \Cal{O}^M \to B^{1, \eps}$ is Lipschitz-continuous on the open subset
\begin{align*}
	\Cal{O}^M
		&:= \set {\bgamma \in B^{1, \eps} \colon \abs{\bgamma}_* > M^{-1}, \,
				\norm{\bgamma'}_{L^\iny} < M}, \\
	\abs{\bgamma}_*
		&:= \inf_{\al \ne \al'} \frac{\bgamma(\al) - \bgamma(\al')}{\abs{\al - \al'}}
\end{align*}
for some $M > 0$. A Picard fixed point theorem (Theorem 8.3 of \cite{MB2002}) then assures a local-in-time solution to the ODE,
\begin{align*}
	\diff{\bgamma}{t} = F(\bgamma), \quad
	\bgamma(0) = \bgamma_0 \in \Cal{O}^M,
\end{align*}
with $\bgamma \in C^1([-T, T]; \Cal{O}^M)$ for a $T$ that depends upon $M$.

To adapt the argument in \cite{MB2002} to Type 2 solutions, we decompose $\log \KInfFunc(\x)$ as follows. Let $\varphi \in C_0^\iny(\Pi)$ be a radially symmetric cutoff function supported on $B_{1/4}(0)$ with $\varphi \equiv 1$ on $B_{1/8}(0)$. Then
\begin{align*}
	\log \KInfFunc(\x)
		&= \varphi(\x) \log \abs{\x} + R(\x), \\
	R(\x)
		&:= \varphi(\x) \brac{\log \KInfFunc(\x) - \log \abs{\x}}
			+ (1 - \varphi(\x)) \log \KInfFunc(\x).
\end{align*}
Recall that on $\Pi$, we use coordinates in which $\x = (x_1, x_2)$ with $-1/2 \le x_1 < 1/2$. Because $\varphi(\x) = 0$ for $\abs{x_1} > 1/4$, the function $\varphi(\x) \log \abs{\x}$ is in $C^\iny(\Pi \setminus (0, 0))$. Also, $\log \KInfFunc(\x)$ is harmonic away from the origin, so $R(\x) \in C^\iny(\Pi)$, as follows from \cref{L:loglogDiff}. In particular, $\varphi(\x) \log \abs{\x}$ and $R(\x)$ are well-defined as functions on $\Pi$.

It follows that for each component of $\prt \Omega_{t, 0}$, $F = F_1 + F_2$, where
\begin{align*}
	F_1(\bgamma(\beta))
		&:= \frac{\omega_0}{2 \pi} \int_0^{2 \pi}
				\varphi({\bgamma(\beta) - \bgamma(\al)})
					\log \abs{\bgamma(\beta) - \bgamma(\al)}
				\prt_\al \bgamma(\al) \, d \al, \\
	F_2(\bgamma(\beta))
		&:= \frac{\omega_0}{2 \pi} \int_0^{2 \pi}
				R(\bgamma(\beta) - \bgamma(\al))
				\prt_\al \bgamma(\al) \, d \al.
\end{align*}
Other than the cutoff function, which introduces no real difficulties, $F_1$ is the same expression as in the classical setting and is estimated in $B^{1, \eps}$ in the same manner. We note that applying $d/d \beta$ to $F_1(\bgamma(\beta))$ leads to a singularity in the integrand at $\al = \beta$. The key to estimating $F_1$ is treating $dF_1/d \beta$, beginning in Lemma 8.7 of \cite{MB2002}, as a principal value integral. The situation is no different here than in \cite{MB2002}.

Similarly, for $F_2$, the key is bounding $dF_2/d \beta$ in $C^\eps$. This is much simpler than bounding $d F_1/d \beta$, for we have
\begin{align*}
	 \diff{}{\beta} F_2(\bgamma(\beta))
		&= \frac{\omega_0}{2 \pi} \int_0^{2 \pi}
			\grad R(\bgamma(\beta) - \bgamma(\al)) \cdot \prt_\beta \bgamma(\beta))
				\, \prt_\al \bgamma(\al) \, d \al.
\end{align*}
Then for any $\al$,
\begin{align*}
	&\norm{\grad R(\bgamma(\beta) - \bgamma(\al)) \cdot \prt_\beta \bgamma(\beta))
				\, \prt_\al \bgamma(\al)}_{C^\eps} \\
		&\qquad
		\le \abs{\prt_\al \bgamma(\al)}
			\norm{\grad R}_{C^\eps(\Pi)}
			\norm{\bgamma(\beta) - \bgamma(\al)}_{lip}^\eps
			\norm{\prt_\beta \bgamma(\beta)}_{C^\eps(0, 2 \pi)}.
\end{align*}
But, $\abs{\prt_\al \bgamma(\al)} \le \norm{\bgamma}_{Lip} < M$ and
$\norm{\bgamma(\beta) - \bgamma(\al)}_{lip} = \norm{\bgamma}_{lip} < M$. Hence,
\begin{align*}
	 \norm[\bigg]{\diff{}{\beta} F_2(\bgamma(\beta))}_{C^\eps(0, 2 \pi)}
		&\le C M^2 \abs{\omega_0} \norm{\bgamma}_{C^\eps}.
\end{align*}
We see, then, that the bounds in Lemma 8.10 of \cite{MB2002} hold, and the proof of local-in-time existence is completed as in \cite{MB2002}.

\bigskip\noindent\textbf{Global-in-time $C^{1, \eps}$ solutions}:
The proof of the global existence of a $C^{1, \eps}$ solution to the CDE is the same as in Section 8.3.3 of \cite{MB2002}, except that \cref{C:graduType} is used to obtain $\grad \uu$. By virtue of \cref{P:SymConv}, the estimates differ little from those for classical vortex patches.

This completes the proof for Type 2 solutions. The result for Types 1 and 3 solutions then follows directly, exploiting the lifting of domains described in \cref{S:Lifting}.
\end{proof}

\appendix

\section{Proof of the formula for $\grad u$}\label{A:VelocityGradient}

\noindent 

Before giving the proof of the singular integral operator formula for $\grad u$ of \cref{L:graduTypeI}, let us calculate $\grad K_\iny(\x)$ to obtain the expression for $\beta$. Letting
\begin{align*}
	\xi(\x)
		&= \KInfFunc(\x)^2
		= \sin^2 (\pi x_1) + \sinh^2(\pi x_2),
\end{align*}
we have $\prt_1 \rho(\x) = \pi \sin(2 \pi x_1)$, $\prt_2 \rho(\x) = \pi \sinh(2 \pi x_2)$.
Then from \cref{L:KInfgradperG}, we have
$G(\x) = (2 \pi)^{-1} \log \KInfFunc(\x) = (4 \pi)^{-1} \log \xi(\x)$, so
\begin{align*}
	K_\iny(\x)
		&= \grad^\perp G(\x)
		= \frac{\grad^\perp \xi(\x)}{4 \pi \xi(\x)}
		= \frac{(-\prt_2 \xi(\x), \prt_1 \xi(\x)}{4 \pi \xi(\x)}
		= \frac{(-\sinh(2 \pi x_2), \sin(2 \pi x_1))}{4 \xi(\x)}.
\end{align*}

\begin{remark}
	As in \cref{R:gradKInfSize}, near the origin,
	$\xi(\x) = \KInfFunc(\x)^2 \approx \pi^2 \abs{\x}^2$.
	Hence, $G(\x) \approx (1/4 \pi) \log (\pi^2 \abs{\x}^2) \approx C + (1/2\pi) \log \abs{x}$,
	like the fundamental solution to the Laplacian on $\R^2$.
	Then
	$
		K_\iny(\x)
			\approx 2 \pi \abs{\x}/(4 \xi(\x))
			\approx 2 \pi \abs{\x}/(4 \pi^2 \abs{\x}^2)
			= 1/(2 \pi \abs{x})
	$,
	as it is for the Biot-Savart kernel on $\R^2$.
\end{remark}

Taking another derivative,
\begin{align*}
	\grad K_\iny(\x)
		&=
		\frac{1}{4}
		\begin{pmatrix}
			-\prt_1 \dfrac{\sinh(2 \pi x_2)}{\xi(\x)} &
				-\prt_2 \dfrac{\sinh(2 \pi x_2)}{\xi(\x)} \\
				\\
			\prt_1 \dfrac{\sin(2 \pi x_1)}{\xi(\x)} &
				\prt_2 \dfrac{\sin(2 \pi x_1)}{\xi(\x)}
		\end{pmatrix} \\ \\
		&=
		- \frac{1}{4 \xi(\x)^2}
		\begin{pmatrix}
			-\sinh(2 \pi x_2) \prt_1 \xi(\x) &
				-\sinh(2 \pi x_2) \prt_2 \xi(\x) \\
				\\
			\sin(2 \pi x_1) \prt_1 \xi(\x) &
				\sin(2 \pi x_1) \prt_2 \xi(\x)
		\end{pmatrix} \\ \\
		&\qquad
		+
		\frac{1}{4 \xi(\x)}
		\begin{pmatrix}
			0 &
				-2 \pi \cosh(2 \pi x_2) \\
				\\
			2 \pi \cos(2 \pi x_1) &
				0
		\end{pmatrix} \\
		&=
		- \frac{1}{4 \xi(\x)^2}
		\begin{pmatrix}
			-\sinh(2 \pi x_2) \pi \sin(2 \pi x_1) &
				-\sinh(2 \pi x_2) \pi \sinh(2 \pi x_2) \\
				\\
			\sin(2 \pi x_1) \pi \sin(2 \pi x_1) &
				\sin(2 \pi x_1) \pi \sinh(2 \pi x_2)
		\end{pmatrix} \\ \\
		&\qquad
		+
		\frac{1}{4 \xi(\x)^2}
		\begin{pmatrix}
			0 &
				-2 \pi \cosh(2 \pi x_2) \xi(\x) \\
				\\
			2 \pi \cos(2 \pi x_1) \xi(\x) &
				0
		\end{pmatrix} \\
		&=
		\frac{\pi}{4 \xi(\x)^2}
		\begin{pmatrix}
			\sinh(2 \pi x_2) \sin(2 \pi x_1) &
				\sinh^2(2 \pi x_2) \\
				\\
			-\sin^2(2 \pi x_1) &
				-\sin(2 \pi x_1) \sinh(2 \pi x_2)
		\end{pmatrix} \\ \\
		&\qquad
		+
		\frac{\pi}{4 \xi(\x)^2}
		\begin{pmatrix}
			0 &
				-2 \cosh(2 \pi x_2) \xi(\x) \\
				\\
			2 \cos(2 \pi x_1) \xi(\x) &
				0
		\end{pmatrix} \\
	&= \frac{\pi}{2 \rho(\x)^4}
		\begin{pmatrix}
			\al_{11}(\x) & \al_{12}(\x) \\
			\al_{21}(\x) & \al_{22}(\x)
		\end{pmatrix},
\end{align*}
where
\begin{align*}
	\al_{11}(\x)
		&= - \al_{22}(\x)
		= \tfrac{1}{2} \sinh(2 \pi x_2) \sin(2 \pi x_1), \\
	\al_{12}(\x)
		&= \tfrac{1}{2}
			\brac{
				\sinh^2(2 \pi x_2)
				- 2 \cosh(2 \pi x_2) \xi(\x)
			}, \\
	\al_{21}(\x)
		&= \tfrac{1}{2}
			\brac{
				-\sin^2(2 \pi x_1)
				+ 2 \cos(2 \pi x_1) \xi(\x)
			}.
\end{align*}
Using \cref{e:TrigIdentities} and $\cosh^2 x - \sinh^2 x = 1$, we see that
\begin{align*}
	2 \al_{12}(\x)
		&= \sinh^2(2 \pi x_2)
				- 2 \cosh(2 \pi x_2) (\sin^2(\pi x_1) + \sinh^2(\pi x_2)) \\
		&= \sinh^2(2 \pi x_2)
				- 2 \cosh(2 \pi x_2) \sin^2(\pi x_1)
				- \cosh(2 \pi x_2) (\cosh(2 \pi x_2) - 1) \\
		&= - 1 + \cosh(2 \pi x_2)(1 - 2\sin^2(\pi x_1))
		= \cosh(2 \pi x_2) \cos(2 \pi x_1) - 1, \\
	2 \al_{21}(\x)
		&= -\sin^2(2 \pi x_1)
				+ 2 \cos(2 \pi x_1) (\sin^2(\pi x_1) + \sinh^2(\pi x_2)) \\
		&= -\sin^2(2 \pi x_1)
				+ \cos(2 \pi x_1) (1 - \cos(2 \pi x_1))
					+ 2 \cos(2 \pi x_1) \sinh^2(\pi x_2)) \\
		&= -1 + \cos(2 \pi x_1) (1 + 2 \sinh^2(\pi x_2))
		= \cos(2 \pi x_1) \cosh(2 \pi x_2) - 1.
\end{align*}

Thus,
\begin{align*}
	\grad K_\iny(\x)
		&= \frac{\pi}{2} \frac{\beta(\x)}{\rho(\x)^2},
\end{align*}
where
\begin{align*}
	\beta(\x)
		&= \frac{1}{2 \rho(\x)^2}
			\begin{pmatrix}
				\sin(2 \pi x_1) \sinh(2 \pi x_2) &
					\cos(2 \pi x_1) \cosh(2 \pi x_2) - 1 \\
				\cos(2 \pi x_1) \cosh(2 \pi x_2) - 1 &
					- \sin(2 \pi x_1) \sinh(2 \pi x_2)
			\end{pmatrix},
\end{align*}
as given in \cref{L:graduTypeI}.

\begin{proof}[\textbf{Proof of \cref{L:graduTypeI}}]
	Let $M \in (H^1(\Omega)^{2 \times 2}$ be arbitrary. We will show that
	\begin{align*}
		(\grad \uu, M)
			&= \pr{\sum_{n \in \Z} \frac{\omega(\x + (n, 0))}{2}
				\begin{pmatrix}
					0 & -1 \\
					1 & 0
				\end{pmatrix}
				, M}
				+ \frac{1}{2 \pi}
				\PV \int_{\R^2} \grad K_\iny(\x - \y) M(\y) \, d \y,
	\end{align*}
	giving the action of $\grad \uu \in H^{-1}(\R^2)$ on any test function
	in $H^1(\R^2)$, and thus establishing our expression for $\grad \uu$.

	For any $r \in (0, 1)$, we let
	\begin{align*}
		U_r
			= \bigcup_{n \in \Z} B_r(\x + (n, 0)).
	\end{align*}

	Then
	\begin{align*}
		(\grad \uu, &M)
			= (\uu, \dv M)
			= (K_\iny * \omega, \dv M) 
			= \lim_{r \to 0} \int_{U_r^C} K_\iny * \omega(\x) \,
					\dv M(\x) \, d \x
				\\
			&= - \lim_{r \to 0} \int_{U_r^C}
					\grad (K_\iny * \omega)(\x) \, M(\x) \, d \x
						- \lim_{r \to 0} \int_{\prt U_r} (\grad M \cdot \n)
							K_\iny * \omega \, d S
			=: I + II.
	\end{align*}
	We used here that $\uu$ is integrable and that the orientation of $\prt U$
	is opposite that of $\prt U^C$. The limit in $I$ gives the principal value integral
	in our expression for $\grad u$.
	Noting that the compact support of $\omega$ makes the sum below finite,
	\begin{align*}
		II
			&= \sum_{n \in \Z}
				\lim_{r \to 0} \int_{\prt B_r(\x + (n, 0))} (\grad M \cdot \n)
							K_\iny * \omega \, d S \\
			&= \sum_{n \in \Z}
				\lim_{r \to 0} \int_{\prt B_r(\x)} (\grad M(\cdot + (n, 0)) \cdot \n)
							K_\iny * \omega \, d S \\
			&= \sum_{n \in \Z}
				\lim_{r \to 0} \int_{\prt B_r(\x)} (\grad M(\cdot + (n, 0)) \cdot \n)
							K * \omega \, d S \\
			&= \sum_{n \in \Z}
				\pr{\frac{\omega}{2}
				\begin{pmatrix}
					0 & -1 \\
					1 & 0
				\end{pmatrix},
				M(\cdot + (n, 0))
				}
			= \sum_{n \in \Z}
				\pr{\frac{\omega(\x - (n, 0))}{2}
				\begin{pmatrix}
					0 & -1 \\
					1 & 0
				\end{pmatrix},
				M
				}.
	\end{align*}
	We used that $K_\iny(\y)$ becomes $K(\y)$ in the limit of small $\y$, and then evaluated
	the limit of the boundary integral as in the classical case.
\end{proof}

\section*{Acknowledgements}
\noindent
DMA is grateful to the National Science Foundation for support through grant DMS-1907684.


\begin{thebibliography}{10}

\bibitem{AbramowitzStegun1964}
Milton Abramowitz and Irene~A. Stegun.
\newblock {\em Handbook of mathematical functions with formulas, graphs, and
  mathematical tables}.
\newblock National Bureau of Standards Applied Mathematics Series, No. 55. U.
  S. Government Printing Office, Washington, D.C., 1964.
\newblock For sale by the Superintendent of Documents.

\bibitem{AfendikovMielke2005}
Andrei~L. Afendikov and Alexander Mielke.
\newblock Dynamical properties of spatially non-decaying 2{D} {N}avier-{S}tokes
  flows with {K}olmogorov forcing in an infinite strip.
\newblock {\em J. Math. Fluid Mech.}, 7(suppl. 1):S51--S67, 2005.

\bibitem{Ahlfors}
Lars~V. Ahlfors.
\newblock {\em Complex analysis}.
\newblock International Series in Pure and Applied Mathematics. McGraw-Hill
  Book Co., New York, third edition, 1978.
\newblock An introduction to the theory of analytic functions of one complex
  variable.

\bibitem{AKLL2015}
David~M. Ambrose, James~P. Kelliher, Milton C.~Lopes Filho, and Helena
  J.~Nussenzveig Lopes.
\newblock {S}erfati solutions to the 2{D} {E}uler equations on exterior
  domains.
\newblock {\em Journal of Differential Equations}, 259(9):4509--4560, 2015.

\bibitem{bernoff}
Oliver~V. Atassi, Andrew~J. Bernoff, and Seth Lichter.
\newblock The interaction of a point vortex with a wall-bounded vortex layer.
\newblock {\em J. Fluid Mech.}, 343:169--195, 1997.

\bibitem{baeKelliher}
Hantaek Bae and James~P. Kelliher.
\newblock Propagation of regularity of level sets for a class of active
  transport equations.
\newblock {\em J. Math. Anal. Appl.}, 497(1):Paper No. 124823, 37, 2021.

\bibitem{benedettoPulvirenti}
D.~Benedetto and M.~Pulvirenti.
\newblock From vortex layers to vortex sheets.
\newblock {\em SIAM J. Appl. Math.}, 52(4):1041--1056, 1992.

\bibitem{bertozziEtAl}
A.~Bertozzi, J.~Garnett, T.~Laurent, and J.~Verdera.
\newblock The regularity of the boundary of a multidimensional aggregation
  patch.
\newblock {\em SIAM J. Math. Anal.}, 48(6):3789--3819, 2016.

\bibitem{ConstantinBertozzi1993}
A.~L. Bertozzi and P.~Constantin.
\newblock Global regularity for vortex patches.
\newblock {\em Comm. Math. Phys.}, 152(1):19--28, 1993.

\bibitem{BertozziThesis1991}
Andrea~Louise Bertozzi.
\newblock {\em Existence, uniqueness, and a characterization of solutions to
  the contour dynamics equation}.
\newblock ProQuest LLC, Ann Arbor, MI, 1991.
\newblock Thesis (Ph.D.)--Princeton University.

\bibitem{caflisch2}
R.~E. Caflisch, F.~Gargano, M.~Sammartino, and V.~Sciacca.
\newblock Complex singularity analysis for vortex layer flows.
\newblock {\em J. Fluid Mech.}, 932:Paper No. A21, 37, 2022.

\bibitem{caflisch1}
R.~E. Caflisch, M.~C. Lombardo, and M.~M.~L. Sammartino.
\newblock Vortex layers of small thickness.
\newblock {\em Comm. Pure Appl. Math.}, 73(10):2104--2179, 2020.

\bibitem{caflisch3}
Russel~E. Caflisch, Francesco Gargano, Marco Sammartino, and Vincenzo Sciacca.
\newblock Complex singularities and {PDE}s.
\newblock {\em Riv. Math. Univ. Parma (N.S.)}, 6(1):69--133, 2015.

\bibitem{chaeEtal}
Dongho Chae, Peter Constantin, Diego C\'{o}rdoba, Francisco Gancedo, and
  Jiahong Wu.
\newblock Generalized surface quasi-geostrophic equations with singular
  velocities.
\newblock {\em Comm. Pure Appl. Math.}, 65(8):1037--1066, 2012.

\bibitem{cheminCRAS}
Jean-Yves Chemin.
\newblock Existence globale pour le probl\`eme des poches de tourbillon.
\newblock {\em C. R. Acad. Sci. Paris S\'{e}r. I Math.}, 312(11):803--806,
  1991.

\bibitem{cheminENS}
Jean-Yves Chemin.
\newblock Persistance de structures g\'{e}om\'{e}triques dans les fluides
  incompressibles bidimensionnels.
\newblock {\em Ann. Sci. \'{E}cole Norm. Sup. (4)}, 26(4):517--542, 1993.

\bibitem{Co1978}
John~B. Conway.
\newblock {\em Functions of one complex variable}, volume~11 of {\em Graduate
  Texts in Mathematics}.
\newblock Springer-Verlag, New York, second edition, 1978.

\bibitem{crowdy}
Darren Crowdy.
\newblock Exact solutions for uniform vortex layers attached to corners and
  wedges.
\newblock {\em European J. Appl. Math.}, 15(6):643--650, 2004.

\bibitem{Gallay2017}
Thierry Gallay.
\newblock Infinite energy solutions of the two-dimensional {N}avier-{S}tokes
  equations.
\newblock {\em Ann. Fac. Sci. Toulouse Math. (6)}, 26(4):979--1027, 2017.

\bibitem{GallaySlijepcevic2014}
Thierry Gallay and Sini\v{s}a Slijep\v{c}evi\'{c}.
\newblock Energy bounds for the two-dimensional {N}avier-{S}tokes equations in
  an infinite cylinder.
\newblock {\em Comm. Partial Differential Equations}, 39(9):1741--1769, 2014.

\bibitem{GallaySlijepcevic2015}
Thierry Gallay and Sini\v{s}a Slijep\v{c}evi\'{c}.
\newblock Uniform boundedness and long-time asymptotics for the two-dimensional
  {N}avier-{S}tokes equations in an infinite cylinder.
\newblock {\em J. Math. Fluid Mech.}, 17(1):23--46, 2015.

\bibitem{gancedo}
Francisco Gancedo.
\newblock Existence for the {$\alpha$}-patch model and the {QG} sharp front in
  {S}obolev spaces.
\newblock {\em Adv. Math.}, 217(6):2569--2598, 2008.

\bibitem{gargano}
Francesco Gargano, Maria~Carmela Lombardo, Marco Sammartino, and Vincenzo
  Sciacca.
\newblock Singularity formation and separation phenomena in boundary layer
  theory.
\newblock In {\em Partial differential equations and fluid mechanics}, volume
  364 of {\em London Math. Soc. Lecture Note Ser.}, pages 81--120. Cambridge
  Univ. Press, Cambridge, 2009.

\bibitem{golubeva}
Natalia~Yurievna Golubeva.
\newblock {\em Singularities in the spatial complex plane for vortex sheets and
  thin vortex layers}.
\newblock ProQuest LLC, Ann Arbor, MI, 2003.
\newblock Thesis (Ph.D.)--The Ohio State University.

\bibitem{hunter1}
John~K. Hunter, Jingyang Shu, and Qingtian Zhang.
\newblock Contour dynamics for surface quasi-geostrophic fronts.
\newblock {\em Nonlinearity}, 33(9):4699--4714, 2020.

\bibitem{hunter2}
John~K. Hunter, Jingyang Shu, and Qingtian Zhang.
\newblock Two-front solutions of the {SQG} equation and its generalizations.
\newblock {\em Commun. Math. Sci.}, 18(6):1685--1741, 2020.

\bibitem{hunter3}
John~K. Hunter, Jingyang Shu, and Qingtian Zhang.
\newblock Global solutions of a surface quasigeostrophic front equation.
\newblock {\em Pure Appl. Anal.}, 3(3):403--472, 2021.

\bibitem{Kelliher2015}
James~P. Kelliher.
\newblock A characterization at infinity of bounded vorticity, bounded velocity
  solutions to the 2{D} {E}uler equations.
\newblock {\em Indiana Univ. Math. J.}, 64(6):1643--1666, 2015.

\bibitem{kiselev}
Alexander Kiselev, Yao Yao, and Andrej Zlato\v{s}.
\newblock Local regularity for the modified {SQG} patch equation.
\newblock {\em Comm. Pure Appl. Math.}, 70(7):1253--1315, 2017.

\bibitem{MB2002}
Andrew~J. Majda and Andrea~L. Bertozzi.
\newblock {\em Vorticity and incompressible flow}, volume~27 of {\em Cambridge
  Texts in Applied Mathematics}.
\newblock Cambridge University Press, Cambridge, 2002.

\bibitem{MP1984}
C.~Marchioro and M.~Pulvirenti.
\newblock {\em Vortex methods in two-dimensional fluid dynamics}, volume 203 of
  {\em Lecture Notes in Physics}.
\newblock Springer-Verlag, Berlin, 1984.

\bibitem{MP1994}
Carlo Marchioro and Mario Pulvirenti.
\newblock {\em Mathematical theory of incompressible nonviscous fluids},
  volume~96 of {\em Applied Mathematical Sciences}.
\newblock Springer-Verlag, New York, 1994.

\bibitem{pullin}
D.I. Pullin and P.A. Jacobs.
\newblock Inviscid evolution of stretched vortex arrays.
\newblock {\em J. Fluid Mech.}, 171:377–406, 1986.

\bibitem{rodrigo1}
Jos\'{e}~Luis Rodrigo.
\newblock The vortex patch problem for the surface quasi-geostrophic equation.
\newblock {\em Proc. Natl. Acad. Sci. USA}, 101(9):2684--2686, 2004.

\bibitem{rodrigo2}
Jos\'{e}~Luis Rodrigo.
\newblock On the evolution of sharp fronts for the quasi-geostrophic equation.
\newblock {\em Comm. Pure Appl. Math.}, 58(6):821--866, 2005.

\bibitem{serfatiStratifee}
Philippe Serfati.
\newblock R\'{e}gularit\'{e} stratifi\'{e}e et \'{e}quation d'{E}uler {$3$}{D}
  \`a temps grand.
\newblock {\em C. R. Acad. Sci. Paris S\'{e}r. I Math.}, 318(10):925--928,
  1994.

\bibitem{Serfati1995Bounded}
Philippe Serfati.
\newblock Solutions {$C^\infty$} en temps, {$n$}-{$\log$} {L}ipschitz
  born\'{e}es en espace et \'{e}quation d'{E}uler.
\newblock {\em C. R. Acad. Sci. Paris S\'{e}r. I Math.}, 320(5):555--558, 1995.

\bibitem{TaniuchiEtAl2010}
Yasushi Taniuchi, Tomoya Tashiro, and Tsuyoshi Yoneda.
\newblock On the two-dimensional {E}uler equations with spatially almost
  periodic initial data.
\newblock {\em J. Math. Fluid Mech.}, 12(4):594--612, 2010.

\end{thebibliography}

\def\cprime{$'$} \def\polhk#1{\setbox0=\hbox{#1}{\ooalign{\hidewidth
  \lower1.5ex\hbox{`}\hidewidth\crcr\unhbox0}}}

\end{document}